\documentclass[12pt, reqno]{amsart}
\usepackage{amsmath, amsthm, amscd, amsfonts, amssymb,mathrsfs, graphicx,enumerate,comment,array}
\usepackage{color}
\usepackage{float}
\usepackage[bookmarksnumbered, colorlinks, plainpages]{hyperref}
\hypersetup{colorlinks=true,linkcolor=red, anchorcolor=green, citecolor=cyan, urlcolor=red, filecolor=magenta, pdftoolbar=true}

\allowdisplaybreaks[4]

\newtheorem{theorem}{Theorem}[section]
\newtheorem{lemma}[theorem]{Lemma}
\newtheorem{proposition}[theorem]{Proposition}
\newtheorem{corollary}[theorem]{Corollary}
\theoremstyle{definition}
\newtheorem{definition}[theorem]{Definition}

\theoremstyle{remark}
\newtheorem{remark}[theorem]{Remark}
\numberwithin{equation}{section}

\begin{document}
\setcounter{page}{1}

\title[Schr\"{o}dinger maximal estimate]{$L^2$ Schr\"{o}dinger maximal estimates associated with finite type phases in $\mathbb{R}^2$}
\author[Z. Li, J. Zhao \MakeLowercase{and} T. Zhao]{Zhuoran Li$^1$, Junyan
Zhao$
^2$$^*$ \MakeLowercase{and} Tengfei Zhao$^3$}
\address{$^1$
The Graduate School of China Academy of Engineering Physics, P. O. Box 2101, Beijing
100088, P.R. China.}
\email{\textcolor[rgb]{0.00,0.00,0.84}{lizhuoran18@gscaep.ac.cn}}
\address{$^2$College of Mathematics and Computer Science, Zhejiang Normal
University, Jinhua 321004, P.R. China.}
\email{\textcolor[rgb]{0.00,0.00,0.84}{jyzhaomath@zjnu.edu.cn}}
\address{$^3$School of Mathematics and Physics,
University of Science and Technology Beijing,
Beijing 100083, P.R. China.}
\email{\textcolor[rgb]{0.00,0.00,0.84}{zhao\_tengfei@ustb.edu.cn}}
\date{Received: xxxxxx; Revised: yyyyyy; Accepted: zzzzzz. \\
\indent $^{*}$Corresponding author}

\begin{abstract}
In this paper, we establish Schr\"{o}dinger maximal estimates associated
with the finite type phases
\begin{equation*}
\phi(\xi_1,\xi_2):=\xi^m_1+\xi^m_2,\;(\xi_1,\xi_2)\in [0,1]^2,
\end{equation*}
where $m \geq 4$ is an even number. Following~\cite{DuZhang19-Annals}, we
prove an $L^2$ fractal restriction estimate associated with the surfaces
\begin{equation*}
F^2_m:=\{(\xi_1,\xi_2,\phi(\xi_1,\xi_2)):\;(\xi_1,\xi_2)\in [0,1]^2\}
\end{equation*}
as the main result, which also gives results on the average Fourier decay of
fractal measures associated with these surfaces. The key ingredients of the
proof include the rescaling technique from \cite{LiMiaoZheng2020},
Bourgain-Demeter's $\ell^2$ decoupling inequality, the reduction of
dimension arguments from \cite{LiZheng21} and induction on scales.
\end{abstract}
\maketitle

\section{\textbf{Introduction}}\label{sect:introd}
The solution to the Cauchy problem of the free Schr\"{o}dinger equation
\begin{equation*}
\left\{
\begin{array}{l}
i\partial _{t}u-\Delta u=0,\,(x,t)\in \mathbb{R}^{n}\times \mathbb{R}\,, \\
u(x,0)=f(x) \\
\end{array}
\right.
\end{equation*}
is given by
\begin{equation*}
e^{it\Delta}f(x)=(2\pi)^{-n}\int_{\mathbb{R}^n} \hat{f}(\xi)e(x\cdot\xi+t
\vert \xi \vert^2)d\xi,
\end{equation*}
where $e(t):=e^{it}$.

Carleson in \cite{Carleson} proposed the problem of determining the optimal $
s$ for which
\begin{equation*}
\lim_{t\rightarrow 0}e^{it\Delta }f(x)=f(x)
\end{equation*}
almost everywhere, whenever $f\in H^{s}(\mathbb{R}^{n})$. The problem has
attracted many authors. In particular, Lee \cite{Lee06} proved that almost
everywhere convergence holds if $s>\frac{3}{8}$ when $n=2$. For general $n$,
Bourgain \cite{Bo16} gave surprising counterexamples showing that the
convergence fails if $s<\frac{n}{2(n+1)}$. Later, Du, Guth and Li \cite
{DuGuthLi} proved that the convergence holds for $s>\frac{1}{3}$ when $n=2$,
which combining with Bourgain's necessary condition solves the Carleson's
problem when $n=2,$ except at the endpoint $s=\frac{1}{3}$. In higher
dimensions ($n\geq 3$), Du and Zhang \cite{DuZhang19-Annals} proved that
almost everywhere convergence holds in the sharp range $s>\frac{n}{2(n+1)}$.
Thus, Carleson's problem is completely solved for any dimension $n,$ except
at the endpoint $\frac{n}{2(n+1)}$. For the study of pointwise convergence
problem on 2D fractional order Schr\"{o}dinger operators, we refer to Miao,
Yang and Zheng \cite{MYZ}. See also C. Cho and H. Ko \cite{ChoKo}.

Let $m \geq 4$ be an even number. The solution to the Cauchy problem of the following generalized
Schr\"{o}dinger equation
\begin{equation*}
\left\{
\begin{array}{l}
i\partial _{t}u-(\partial^m_1u+\partial^m_2u)=0,\,(x,t)\in \mathbb{R}
^{2}\times \mathbb{R}\,, \\
u(x,0)=f(x) \\
\end{array}
\right.
\end{equation*}
is given by
\begin{equation*}
e^{it\phi(D)}f(x)=(2\pi)^{-2}\int_{\mathbb{R}^2} \hat{f}(\xi)e[x_1\xi_1+x_2
\xi_2+t\phi(\xi)]d\xi,
\end{equation*}
where $\phi(\xi_1,\xi_2):=\xi^m_1+\xi^m_2$.

We note that the Gaussian curvature of the surfaces associated with $
\phi(\xi)$ vanishes when $\xi_1=0$ or $\xi_2=0$, which is different from
those non-degenerate cases in the literature. To our knowledge, there are
few works on the sufficient conditions for the almost everywhere convergence
problem associated with such degenerate phases. For the study of necessary
conditions, we refer to \cite{AnChuPierce21,EV2021}.

In this article, we establish the following result.

\begin{theorem}
\label{thm:main1} For every $f\in H^s(\mathbb{R}^2)$ with $s>\frac{1}{2}-
\frac{1}{3m}$, $\displaystyle\lim_{t\rightarrow 0}e^{it\phi(D)}f(x)=f(x)$
almost everywhere.
\end{theorem}

We use $B^{d}(x,r)$ to denote the ball centered at $x$ with radius $r$ in $
\mathbb{R}^{d}$. By a standard approximation argument, Theorem \ref
{thm:main1} is a consequence of the following Schr\"{o}dinger maximal
estimate associated with $\phi $:

\begin{theorem}
\label{thm:main2} For any $s>\frac{1}{2}-\frac{1}{3m}$, there holds
\begin{equation*}
\Vert \sup_{0<t \leq 1}\vert e^{it\phi(D)}f \vert \Vert_{L^2(B^2(0,1))}\leq
C_{s}\Vert f \Vert_{H^s(\mathbb{R}^2)}
\end{equation*}
for any function $f\in H^s(\mathbb{R}^2)$.
\end{theorem}

By Lemma 2.1 in \cite{CLV}, the Littlewood-Paley decomposition and a
rescaling argument, Theorem \ref{thm:main2} is reduced to the following
result.

\begin{theorem}
\label{reducedmaximal} For any $\varepsilon >0$, there exists a constant $
C_{\varepsilon }$ such that
\begin{equation*}
\Vert \sup_{0<t\leq R}|e^{it\phi (D)}f|\Vert _{L^{2}(B^{2}(0,R))}\leq
C_{\varepsilon }R^{\frac{1}{2}-\frac{1}{3m}+\varepsilon }\Vert f\Vert _{2}
\end{equation*}
holds for all $R\geq 1$ and all $f$ with $supp\hat{f}\subset
A(1):=[0,1]^{2}-[0,\frac{1}{2}]^{2}$.
\end{theorem}

\begin{remark}
\label{rek:oddcase} If $m$ is odd, we can not reduce $supp\hat{f}\subset
[-1,1]^2-[-\frac{1}{2},\frac{1}{2}]^2$ to $supp\hat{f}\subset [0, 1]^2-[0,
\frac{1}{2}]^2$. To see it, let us consider the case $supp\hat{f}\subset [
\frac{1}{2},1]\times [-1,-\frac{1}{2}]=:A(1^{\prime })$. Taking the change
of variables
\begin{equation*}
\xi_1=\eta_1, \xi_2=-\eta_2,
\end{equation*}
the phase function $\phi(\eta_1,-\eta_2):=\eta^m_1-\eta^m_2$ can be viewed as
small perturbations of $\phi_0(\eta_1,\eta_2):=\eta^2_1-\eta^2_2$, where
$(\eta_1,\eta_2)\in[0,1]^2$. But we know from \cite{RVV} that
for any $f\in H^s(\mathbb{R}^2)$
\begin{equation*}
\lim_{t\rightarrow 0}e^{it\phi_0(D)}f(x)=f(x),\quad a.e. \ x\in \mathbb{R}^2
\end{equation*}
fails for $s<\frac{1}{2}$. For this reason we only consider that $m$ is even
in the current paper.
\end{remark}

Cubes of the form $m+[0,M]^{3}$ with $m\in (M\mathbb{Z})^{3}$ are called
lattice $M-$cubes. We will prove the following $L^{2}$ restriction estimate,
from which Theorem \ref{reducedmaximal} follows.

\begin{theorem}
\label{fractalrestriction} Suppose that $X=\bigcup_{k}B_{k}$ is a union of
lattice unit cubes in $Q_{R}:=[0,R]^{3}$ and each $R^{\frac{1}{2}}$-cube
intersecting $X$ contains $\sim \nu $ many unit cubes in $X$. Let $1\leq
\alpha \leq 3$ and $\lambda $ be
\begin{equation*}
\lambda :=\max_{B^{3}(x^{\prime },r)\subset Q_{R},x^{\prime }\in \mathbb{R}
^{3},r\geq 1}\frac{\sharp \{B_{k}:B_{k}\subset B^{3}(x^{\prime },r)\}}{
r^{\alpha }}.
\end{equation*}
For any $\varepsilon >0$, there exists a positive constant $C_{\varepsilon }$
such that the following inequalities hold for all $R\geq 1:$
\begin{equation}
\Vert e^{it\phi (D)}f\Vert _{L^{2}(X)}\leq C_{\varepsilon }\lambda ^{\frac{1
}{6}}\nu ^{\frac{1}{6}}R^{\frac{1}{3}-\frac{4-\alpha }{6m}+\varepsilon
}\Vert f\Vert _{2}  \label{a1estimate}
\end{equation}
for all $f$ with $supp\hat{f}\subset A(1)$, and
\begin{equation}
\Vert e^{it\phi (D)}f\Vert _{L^{2}(X)}\leq C_{\varepsilon }\lambda ^{\frac{1
}{6}}\nu ^{\frac{1}{6}}R^{\frac{5}{12}-\frac{5-\alpha }{6m}+\varepsilon
}\Vert f\Vert _{2}  \label{b1estimate}
\end{equation}
for all $f$ with $supp\hat{f}\subset \lbrack 0,1]^{2}$.
\end{theorem}

\begin{remark}
To deduce Theorem \ref{reducedmaximal}, \eqref{a1estimate} is sufficient. We
only use \eqref{b1estimate} in Section 4. The two different behaviors in
\eqref{a1estimate} and \eqref{b1estimate} depend on whether the support of $
\hat{f}$ contains the origin or not. When the support of $\hat{f}$ does not
contain the origin, at most one principal curvature of the surface vanishes.
If the support of $\hat{f}$ contains the origin, both principal curvatures
of the surface studied vanish. Therefore, we deduce a better result in
\eqref{a1estimate} than in \eqref{b1estimate}.
\end{remark}

Clearly, one has $\nu \leq \lambda R^{\alpha/2}$ in Theorem \ref
{fractalrestriction}. As a direct result of Theorem \ref{fractalrestriction},
there holds a weaker $L^2$ restriction estimate:

\begin{corollary}
\label{weakfractal} Suppose that $X=\bigcup_{k}B_{k}$ is a union of lattice
unit cubes in $Q_{R}$. Let $1\leq \alpha \leq 3$ and $\lambda $ be
\begin{equation*}
\lambda :=\max_{B^{3}(x^{\prime },r)\subset Q_{R},x^{\prime }\in \mathbb{R}
^{3},r\geq 1}\frac{\sharp \{B_{k}:B_{k}\subset B^{3}(x^{\prime },r)\}}{
r^{\alpha }}.
\end{equation*}
For any $\varepsilon >0$, there exists a positive constant $C_{\varepsilon }$
such that the following inequalities hold for all $R\geq 1$
\begin{equation*}
\Vert e^{it\phi (D)}f\Vert _{L^{2}(X)}\leq C_{\varepsilon }\lambda ^{\frac{1
}{3}}R^{\frac{1}{3}+\frac{\alpha }{12}-\frac{4-\alpha }{6m}+\varepsilon
}\Vert f\Vert _{2}
\end{equation*}
for all $f$ with $supp\hat{f}\subset A(1)$, and
\begin{equation*}
\Vert e^{it\phi (D)}f\Vert _{L^{2}(X)}\leq C_{\varepsilon }\lambda ^{\frac{1
}{3}}R^{\frac{5}{12}+\frac{\alpha }{12}-\frac{5-\alpha }{6m}+\varepsilon
}\Vert f\Vert _{2}
\end{equation*}
for all $f$ with $supp\hat{f}\subset \lbrack 0,1]^{2}$.
\end{corollary}

Corollary \ref{weakfractal} is sufficient to derive the $L^{2}$ Schr\"{o}
dinger maximal estimate (Theorem \ref{reducedmaximal}). For the details, we
refer to the proof of Theorem 1.3 by Corollary 1.7 in Du-Zhang \cite
{DuZhang19-Annals}. Theorem \ref{fractalrestriction} is based on Proposition
\ref{mainprop} below, which will be proved by an induction.

A collection of quantities are said to be essentially constant provided that
all the quantities lie in the same interval of the form $[2^n,2^{n+1}]$,
where $n\in \mathbb{Z}$.

\begin{proposition}
\label{mainprop} Suppose that $Y=\bigcup_{k=1}^{N}B_{k}$ is a union of
lattice $K$-cubes in $Q_{R}$ and each $R^{\frac{1}{2}}$-cube intersecting $Y$
contains $\sim \nu $ many $K$-cubes in $Y$, where $K=R^{\delta }$. Suppose
that $\Vert e^{it\phi (D)}f\Vert _{L^{6}(B_{k})}$ is essentially constant in
$k=1,2,...,N$. Let $1\leq \alpha \leq 3$ and $\lambda $ be
\begin{equation*}
\lambda :=\max_{B^{3}(x^{\prime },r)\subset Q_{R},x^{\prime }\in \mathbb{R}
^{3},r\geq K}\frac{\sharp \{B_{k}:B_{k}\subset B^{3}(x^{\prime },r)\}}{
r^{\alpha }}.
\end{equation*}
For any $\varepsilon >0$, there exists a positive constant $C_{\varepsilon }$
and $\delta =\varepsilon ^{100}$ such that the following inequalities hold
for all $R\geq 1:$
\begin{equation*}
\Vert e^{it\phi (D)}f\Vert _{L^{6}(Y)}\leq C_{\varepsilon }N^{-\frac{1}{3}
}\lambda ^{\frac{1}{6}}\nu ^{\frac{1}{6}}R^{\frac{1}{3}-\frac{4-\alpha }{6m}
+\varepsilon }\Vert f\Vert _{2}
\end{equation*}
for all $f$ with $supp\hat{f}\subset A(1)$, and
\begin{equation*}
\Vert e^{it\phi (D)}f\Vert _{L^{6}(Y)}\leq C_{\varepsilon }N^{-\frac{1}{3}
}\lambda ^{\frac{1}{6}}\nu ^{\frac{1}{6}}R^{\frac{5}{12}-\frac{5-\alpha }{6m}
+\varepsilon }\Vert f\Vert _{2}
\end{equation*}
for all $f$ with $supp\hat{f}\subset \lbrack 0,1]^{2}$.
\end{proposition}

\begin{remark}
\label{Q1} In the proof of Proposition \ref{mainprop}, we will use the
$\ell^2$ decoupling inequalities for the perturbed parabola and for the curve
\begin{equation*}
\{(s,s^m):s\in [0,1]\}.
\end{equation*}
This is already known due to \cite{Xi, Yang21}. To establish Proposition \ref
{mainprop} in higher dimensions ($n\geq 3$) with the methods of Du-Zhang
\cite{DuZhang19-Annals} and the current paper, one needs the $\ell^2$
decoupling inequalities for the hypersurfaces
\begin{equation*}
\Big\{(\xi_1,...,\xi_{n-1},\phi_1(\xi_1)+...+\phi_s(\xi_s)+\xi_{s+1}^m+...+
\xi_{n-1}^m): (\xi_1,...,\xi_{n-1}) \in [0,1]^{n-1}\Big\},
\end{equation*}
$0 \leq s \leq n-1$, with $\phi_1,...,\phi_s$ being non-degenerate.
With these in hand, one can obtain pointwise convergence results associated
with corresponding finite type phases in higher dimensions. This is
achieved in a recent preprint \cite{GaoLiZhaoZheng2022preprint}.
\end{remark}

From \cite{DuZhang19-Annals}, we know that Proposition \ref{mainprop}
implies Theorem \ref{fractalrestriction}, and therefore it suffices to prove
Proposition \ref{mainprop}. Now we outline the proof of Proposition \ref
{mainprop}.

\vskip 0.2in

\textbf{Outline of the proof of Proposition \ref{mainprop}.} We use
$F_{m}^{2}$ to denote the surface
\begin{equation*}
\{(\xi _{2},\xi _{2},\xi _{1}^{m}+\xi _{2}^{m}):(\xi _{1},\xi _{2})\in
\lbrack 0,1]^{2}\}.
\end{equation*}
The surface $F_{m}^{2}$ is badly behaved when $\xi _{1}=0$ or $\xi _{2}=0$.
The strategy is to single out small neighborhoods of these two lines where
we apply reduction of dimension arguments. More precisely, we first divide
$[0,1]^{2}$ into $\bigcup_{j=0}^{3}\Omega _{j}$, where
\begin{equation*}
\Omega _{0}:=[K^{-\frac{1}{m}},1]\times \lbrack K^{-\frac{1}{m}
},1],\;\;\Omega _{1}:=[K^{-\frac{1}{m}},1]\times \lbrack 0,K^{-\frac{1}{m}}],
\end{equation*}
\begin{equation*}
\Omega _{2}:=[0,K^{-\frac{1}{m}}]\times \lbrack K^{-\frac{1}{m}
},1],\;\;\Omega _{3}:=[0,K^{-\frac{1}{m}}]\times \lbrack 0,K^{-\frac{1}{m}}].
\end{equation*}
Here $K=R^{\delta }$ is a large number. Denote $\hat{f}\mid _{\Omega _{j}}$
by $\hat{f}_{\Omega _{j}}$ for $j=0,1,2,3$. Then, we have
\begin{equation}
\Vert e^{it\phi (D)}f\Vert _{L^{6}(Y)}\leq \sum_{j=0}^{3}\Vert e^{it\phi
(D)}f_{\Omega _{j}}\Vert _{L^{6}(Y)}.  \label{decompositionofthefunction}
\end{equation}
We will estimate $\Vert e^{it\phi (D)}f\Vert _{L^{6}(Y)}$ in several cases.
In \eqref{decompositionofthefunction} one of the
\[\Vert e^{it\phi (D)}f_{\Omega _{j}}\Vert_{L^{6}(Y)}\]
will dominate, and we call it $\Omega _{j}$-case.

We first consider the $\Omega_0$-case. The surface
\begin{equation*}
\{(\xi_2,\xi_2,\xi^m_1+\xi^m_2):(\xi_1,\xi_2)\in \Omega_0\}
\end{equation*}
has two positive principal curvatures with lower bound $K^{-C_0}$ for some
constant $C_0$. For $(\xi_1,\xi_2)\in \Omega_0$, our phase function
$\phi(\xi_1,\xi_2)$ can be viewed as perturbations of $\xi^2_1+\xi^2_2$. By
Proposition 3.1 of Du-Zhang \cite{DuZhang19-Annals}, in the $\Omega_0$-case
we have
\begin{equation*}  \label{omega00}
\Vert e^{it \phi(D)}f \Vert_{L^6(Y)}\leq C_{\varepsilon}K^{O(1)}N^{-\frac{1}{
3}}\lambda^{\frac{1}{6}}\nu^{\frac{1}{6}}R^{\frac{\alpha}{12}
+\varepsilon}\Vert f \Vert_2,
\end{equation*}
where $K^{O(1)}$ is a fixed power of $K$, and the details will be given in
Section 2. For simplicity, we assume that $supp\hat{f}\subset A(1)$. Then,
$\hat{f}_{\Omega_3}=0$. It remains to discuss the $\Omega_1$-case and the
$\Omega_2$-case. By symmetry, it suffices to treat the $\Omega_1$-case. Since
$supp\hat{f}\subset A(1)$, we only need to consider the following subregion
of $\Omega_1$:
\begin{equation*}
\tilde{\Omega}_1:=[\frac{1}{2}, 1]\times [0, K^{-1/m}].
\end{equation*}
We will adapt Bourgain-Demeter's $\ell^2$ decoupling inequality \cite{BD15}
to our needs and reduce the problem to each $K^{-1/2}\times K^{-1/m}$
-rectangle $\tau$. Now it is natural to employ certain rescaling technique.
The main difficulty is that the phase function
\begin{equation*}
\phi(\xi_1,\xi_2):=\xi^m_1+\xi^m_2,\;(\xi_1,\xi_2)\in \tau \subset \tilde{
\Omega}_1
\end{equation*}
is not closed under the change of variables
\begin{equation*}
\left\{
\begin{array}{l}
\xi_1=a+K^{-1/2}\eta_1, \\
\xi_2=K^{-1/m}\eta_2.
\end{array}
\right.
\end{equation*}
In fact, it becomes a new phase function
\begin{equation*}
\phi(a+K^{-1/2}\eta_1, K^{-1/m}\eta_2)\sim \frac{1}{K}(\phi_1(\eta_1)+
\eta_2^m)=\frac{1}{K}\psi(\eta_1,\eta_2),
\end{equation*}
where
\begin{equation*}
\psi(\eta_1,\eta_2):=\phi_1(\eta_1)+\eta_2^m,\;(\eta_1,\eta_2)\in [0,1]^2
\end{equation*}
satisfying
\begin{equation*}
\phi_1^{(\prime\prime)}\sim 1;\quad \vert \phi_1^{(k)}\vert \lesssim 1,
3\leq k \leq m;\quad \phi_1^{(l)}=0, l\geq m+1
\end{equation*}
on the interval $[0,1]$. To overcome the difficulty, we establish a fractal
$L^2$ restriction estimate associated with the new phase function $\psi$ as
an auxiliary proposition (Proposition \ref{mainlemma}), whose proof is given
in Section 3. The key observation is that the new phase function
\begin{equation*}
\psi(\eta_1,\eta_2),\;(\eta_1,\eta_2)\in [0,K^{-1/2}]\times [0, K^{-1/m}]
\end{equation*}
is closed under the change of variable
\begin{equation*}
\left\{
\begin{array}{l}
\xi_1 = K^{-\frac{1}{2}}\eta_1, \\
\xi_2 = K^{-\frac{1}{m}}\eta_2.
\end{array}
\right.
\end{equation*}
With this in hand, we can deduce the desired estimate in the $\Omega_1$
-case. This completes the outline of the proof of Proposition \ref{mainprop}.

\vskip 0.2in

The paper is organized as follows. In Section 2 and Section 3, we give the
proof of Proposition \ref{mainprop} by combining the methods in \cite
{DuZhang19-Annals} and the reduction of dimension arguments from ~\cite
{LiZheng21}. In section 4, we will give two applications of Corollary ~\ref
{weakfractal}.

\vskip 0.2in

\textbf{Notations:} For nonnegative quantities $X$ and $Y$, we will write $
X\lesssim Y$ to denote the estimate $X\leq C Y$ for some large constant $C$
which may vary from line to line and depend on various parameters. If $
X\lesssim Y\lesssim X$, we simply write $X\sim Y$. Dependence of implicit
constants on the power $p$ or the dimension will be suppressed; dependence
on additional parameters will be indicated by subscripts. For example, $
X\lesssim_u Y$ indicates $X\leq CY$ for some $C=C(u)$. For any set $E
\subset \mathbb{R}^d$, we use $\chi_{E}$ to denote the characteristic
function on $E$. Usually, Fourier transform on $\mathbb{R}^d$ is defined by
\begin{equation*}
\widehat{f}(\xi):= (2\pi)^{-d}\int_{\mathbb{R}^d}e^{- ix\cdot \xi}f(x)\,dx.
\end{equation*}

\section{\textbf{The Proof of Proposition \ref{mainprop}}}

Recall that $[0,1]^{2}$ is divided into $\bigcup_{j=0}^{3}\Omega _{j}$,
where
\begin{equation*}
\Omega _{0}:=[K^{-\frac{1}{m}},1]\times \lbrack K^{-\frac{1}{m}},1],\qquad
\Omega _{1}:=[K^{-\frac{1}{m}},1]\times \lbrack 0,K^{-\frac{1}{m}}],
\end{equation*}
\begin{equation*}
\Omega _{2}:=[0,K^{-\frac{1}{m}}]\times \lbrack K^{-\frac{1}{m}},1],\qquad
\Omega _{3}:=[0,K^{-\frac{1}{m}}]\times \lbrack 0,K^{-\frac{1}{m}}].
\end{equation*}
We denote by
\begin{equation*}
\Sigma _{j}:=\{\xi _{1},\xi _{2},\phi (\xi _{1},\xi _{2}):\;(\xi _{1},\xi
_{2})\in \Omega _{j}\}
\end{equation*}
with $j\in \{0,1,2,3\}$. For technical reasons, $K^{-\frac{1}{m}}$, $R^{-\frac{1}{m}}$
and $(\frac{R}{K})^{-1/m}$ should be dyadic numbers. Therefore, we choose $
K=2^{ml}$ and $R=2^{mt}$ $(l,t\in \mathbb{N})$ to be large numbers
satisfying $K\approx R^{\delta }$ with $\delta $ being in Proposition \ref
{mainprop}. Recall that in~\eqref{decompositionofthefunction} we have
\begin{equation*}
\Vert e^{it\phi (D)}f\Vert _{L^{6}(Y)}\leq \sum_{j=0}^{3}\Vert e^{it\phi
(D)}f_{\Omega _{j}}\Vert _{L^{6}(Y)}.
\end{equation*}
We will estimate $\Vert e^{it\phi (D)}f\Vert _{L^{6}(Y)}$ in several cases.
If the contribution from the region $\Omega _{j}$ dominates, we call it
$\Omega_{j}$-case. By the triangle inequality, we have
\begin{equation*}
\Vert e^{it\phi (D)}f_{\Omega _{j}}\Vert _{L^{6}(B)}\geq \frac{1}{16}\Vert
e^{it\phi (D)}f\Vert _{L^{6}(B)}
\end{equation*}
for at least one of the $j$'s.

Now we sort the $K$-cubes in $Y$. Denote
\begin{equation*}
\{B \subset Y: \Vert e^{it\phi(D)}f_{\Omega_j} \Vert_{L^6(B)} \geq \frac{1}{
16}\Vert e^{it\phi(D)}f \Vert_{L^6(B)}\}
\end{equation*}
by $Y^j$. Clearly, one has $Y=\bigcup^{3}_{j=0}Y^j$. If $\sharp \{B: B
\subset Y^j \} \geq \frac{N}{4}$, we call it $\Omega_j$-case.

For the $\Omega_0$-case, we sort the $K$-cubes in $Y^0$ as follows:

\begin{enumerate}
[(1)]
\item
For any dyadic number $A^{(0)}$, let $Y^0_{A^{(0)}}$ be the union of
the $K$-cubes in $Y^0$ satisfying
\begin{equation*}
\Vert e^{it \phi(D)}f_{\Omega_0} \Vert_{L^6(B)}\approx A^{(0)}.
\end{equation*}

\item Fix $A^{(0)}$, for any dyadic number $\nu^{(0)}$, let $
Y^0_{A^{(0)},\nu^{(0)}}$ be the union of the $K$-cubes in $Y^0_{A^{(0)}}$
such that for each $B\subset Y^0_{A^{(0)},\nu^{(0)}}$, the $R^{1/2}$-cube
intersecting $B$ contains $\approx \nu^{(0)}$ cubes from $Y^0_{A^{(0)}}$.
\end{enumerate}

Without loss of generality, we may assume $\Vert f\Vert _{2}=1$. Here $\nu
^{(0)}$ satisfies $1\leq \nu ^{(0)}\lesssim R^{3}$, and the dyadic number $
A^{(0)}$ making significant contributions can be assumed to be between $
R^{-C}$ and $R^{C}$ for a large constant $C$. In fact, if we denote
\begin{equation*}
\{B\subset Y^{0}:\Vert e^{it\phi (D)}f_{\Omega _{0}}\Vert _{L^{6}(B)}\leq
R^{-100}\}
\end{equation*}
by $Y_{-1}^{0}$, then
\begin{equation*}
\Vert e^{it\phi (D)}f_{\Omega _{0}}\Vert _{L^{6}(Y_{-1}^{0})}\lesssim
R^{-90}\Vert f\Vert _{2}.
\end{equation*}
Note that
\begin{equation*}
N^{(0)}\lesssim R^{3},\lambda ^{(0)}\geq R^{-3},\nu ^{(0)}\geq 1,\alpha >0
\end{equation*}
and $K=R^{\delta },\delta =\varepsilon ^{100}$, where
\[N^{(0)}:=\sharp \{B:B\subset Y_{A^{(0)},\nu ^{(0)}}^{0}\}\]
and
\begin{equation*}
\lambda ^{(0)}:=\max_{B^{3}(x^{\prime },r)\subset Q_{R},x^{\prime }\in
\mathbb{R}^{3},r\geq 1}\frac{\sharp \{B_{k}\subset Y^{0}:B_{k}\subset
B^{3}(x^{\prime },r)\}}{r^{\alpha }}.
\end{equation*}
So we have
\begin{equation*}
\Vert e^{it\phi (D)}f_{\Omega _{0}}\Vert _{L^{6}(Y_{-1}^{0})}\leq
C_{\varepsilon }K^{O(1)}(N^{(0)})^{-\frac{1}{3}}(\lambda ^{(0)})^{\frac{1}{6}
}(\nu ^{(0)})^{\frac{1}{6}}R^{\frac{\alpha }{12}+\varepsilon }\Vert f\Vert
_{2}.
\end{equation*}
Therefore, there exist some dyadic numbers $A^{(0)},\nu ^{(0)}$ such that
\begin{equation*}
N^{(0)}\gtrsim \frac{N}{
(\log R)^{2}}.
\end{equation*}
Fix that choice of $A^{(0)},\nu ^{(0)}$ and denote $Y_{A^{(0)},\nu
^{(0)}}^{0}$ by $Y^{0}$ for convenience. Since $\Vert e^{it\phi (D)}f\Vert
_{L^{6}(B_{k})}$ is essentially constant in $k=1,2,...,N$ and $\frac{N}{
(\log R)^{2}}\lesssim N^{(0)}$, we have
\begin{equation}
\Vert e^{it\phi (D)}f\Vert _{L^{6}(Y)}\lesssim (\log R)^{O(1)}\Vert
e^{it\phi (D)}f\Vert _{L^{6}(Y^{0})}.  \label{c18}
\end{equation}
Note that
\begin{equation*}
\Vert e^{it\phi (D)}f_{\Omega _{0}}\Vert _{L^{6}(B)}\geq \frac{1}{16}\Vert
e^{it\phi (D)}f\Vert _{L^{6}(B)}
\end{equation*}
for $B\subset Y^{0}.$ It follows that
\begin{equation*}
\Vert e^{it\phi (D)}f\Vert _{L^{6}(Y^{0})}\lesssim \Vert e^{it\phi
(D)}f_{\Omega _{0}}\Vert _{L^{6}(Y^{0})}.
\end{equation*}
This together with \eqref{c18} gives
\begin{equation}
\Vert e^{it\phi (D)}f\Vert _{L^{6}(Y)}\lesssim _{\varepsilon }R^{\varepsilon
}\Vert e^{it\phi (D)}f_{\Omega _{0}}\Vert _{L^{6}(Y^{0})}  \label{c181}
\end{equation}
in the $\Omega _{0}$-case. We also have
\begin{equation}
\frac{N}{(\log R)^{2}}\lesssim N^{(0)}\leq N,\quad \nu ^{(0)}\leq \nu ,\quad
\lambda ^{(0)}\leq \lambda. \label{parameteromege0}
\end{equation}
Using the argument from \cite{DuZhang19-Annals} as well as the rescaling
trick from Case (c) in Subsection 2.2 of \cite{LiMiaoZheng2020}, we claim
\begin{equation}
\Vert e^{it\phi (D)}f_{\Omega _{0}}\Vert _{L^{6}(Y^{0})}\leq C_{\varepsilon
}K^{O(1)}(N^{(0)})^{-\frac{1}{3}}(\lambda ^{(0)})^{\frac{1}{6}}(\nu ^{(0)})^{
\frac{1}{6}}R^{\frac{\alpha }{12}+\varepsilon }\Vert f\Vert _{2}.
\label{omega0}
\end{equation}
Note that $\min \{\frac{5}{12}-\frac{5-\alpha }{6m},\frac{1}{3}-\frac{
4-\alpha }{6m}\}\geq \frac{\alpha }{12}$ for any $\alpha \in (0,3]$ and each
$m\geq 2$. The factor $K^{O(1)}$ appears on the right-hand side of
\eqref{omega0} because the principal curvature (in the direction of
$\xi_{2} $-axis) of the surface
\begin{equation*}
\Sigma _{0}:=\{(\xi _{1},\xi _{2},\phi (\xi _{1},\xi _{2})):(\xi _{1},\xi
_{2})\in \Omega _{0}\}
\end{equation*}
has a lower bound $K^{-C}$. Inserting ~\eqref{parameteromege0} into ~
\eqref{omega0}, we complete the proof of Proposition \ref{mainprop} in the $
\Omega _{0}$-case.

Now we turn to prove \eqref{omega0}. We divide the region $\Omega _{0}$ into
a family of subregions
\begin{equation*}
\Omega _{0}=\bigcup_{\sigma _{1},\sigma _{2}}\Omega _{\sigma _{1},\sigma
_{2}}
\end{equation*}
for dyadic numbers $\sigma _{i}\in \lbrack K^{-1/m},\frac{1}{2}]$ with
$i\in \{1,2\}$. Here $\Omega _{\sigma _{1},\sigma _{2}}
:=I_{\sigma _{1}}\times I_{\sigma_{2}}$ with $I_{\sigma _{i}}
=[\sigma _{i},2\sigma _{i}]$. We divide
$I_{\sigma _{i}}$ further into
\begin{equation*}
I_{\sigma _{i}}=\bigcup_{j=1}^{\sigma _{i}^{\frac{m}{2}}K^{1/2}}I_{\sigma
_{i},j}.
\end{equation*}
Each $I_{\sigma _{i},j}$ has length of $\sigma _{i}^{-\frac{m-2}{2}}K^{-1/2}$
. So we can write
\begin{equation*}
\Omega _{\sigma _{1},\sigma _{2}}=\bigcup_{\tau }\tau ,
\end{equation*}
where each $\tau $ in the above union has the form
\begin{equation*}
\tau =[a_{1},a_{1}+\sigma _{1}^{-\frac{m-2}{2}}K^{-1/2}]\times \lbrack
a_{2},a_{2}+\sigma _{2}^{-\frac{m-2}{2}}K^{-1/2}]
\end{equation*}
with $\sigma _{i}\leq a_{i}\leq 2\sigma _{i}\;(i=1,2)$. Therefore, we have
\begin{equation}
e^{it\phi (D)}f_{\Omega _{0}}=\sum_{\sigma _{1},\sigma _{2}}e^{it\phi
(D)}f_{\Omega _{\sigma _{1},\sigma _{2}}}, \label{4decompose0}
\end{equation}
and
\begin{equation*}
e^{it\phi (D)}f_{\Omega _{\sigma _{1},\sigma _{2}}}=\sum_{\tau \subset
\Omega _{\sigma _{1},\sigma _{2}}}e^{it\phi (D)}f_{\tau }.
\end{equation*}
We will estimate $\Vert e^{it\phi (D)}f_{\Omega _{0}}\Vert _{L^{6}(Y^{0})}$
in several cases. If the contribution from the region $\Omega _{\sigma
_{1},\sigma _{2}}$ dominates, we call it $\Omega _{\sigma _{1},\sigma _{2}}$
-case. More precisely, given a $K$-cube $B$, by the triangle inequality, we
have
\begin{equation*}
\Vert e^{it\phi (D)}f_{\Omega _{\sigma _{1},\sigma _{2}}}\Vert
_{L^{6}(B)}\gtrsim \frac{\Vert e^{it\phi (D)}f_{\Omega _{0}}\Vert _{L^{6}(B)}
}{(\log K)^{2}}
\end{equation*}
for at least one pair $(\sigma _{1},\sigma _{2})$.

We sort the $K$-cubes in $Y^0$ as follows. Denote by
\begin{equation*}
Y^{\sigma_1,\sigma_2}:=\{B \subset Y^0:\;\Vert e^{it
\phi(D)}f_{\Omega_{\sigma_1,\sigma_2}} \Vert_{L^6(B)}\gtrsim \frac{1}{(\log
K)^2}\Vert e^{it \phi(D)}f_{\Omega_0}\Vert_{L^6(B)}\}.
\end{equation*}
Clearly, one has $Y=\bigcup_{\sigma_1,\sigma_2}Y^{\sigma_1,\sigma_2}$. If
\begin{equation*}
\sharp \{B: B \subset Y^{\sigma_1,\sigma_2}\} \gtrsim \frac{N^{(0)}}{(\log
K)^2},
\end{equation*}
we call it $\Omega_{\sigma_1,\sigma_2}$-case.

For the $\Omega_{\sigma_1,\sigma_2}$-case, we sort the $K$-cubes in $
Y^{\sigma_1,\sigma_2}$ as follows:

\begin{enumerate}
[(1)]
\item For any dyadic number $A^{(\sigma_1,\sigma_2)}$, let $
Y^{\sigma_1,\sigma_2}_{A^{(\sigma_1,\sigma_2)}}$ be the union of the $K$
-cubes in $Y^{\sigma_1,\sigma_2}$ satisfying
\begin{equation*}
\Vert e^{it \phi(D)}f_{\Omega_{\sigma_1,\sigma_2}} \Vert_{L^6(B)}\approx
A^{(\sigma_1,\sigma_2)}.
\end{equation*}

\item Fix $A^{(\sigma_1,\sigma_2)}$, for any dyadic numbers $
\nu^{(\sigma_1,\sigma_2)}$, let $Y^{\sigma_1,\sigma_2}_{A^{(\sigma_1,
\sigma_2)},\nu^{(\sigma_1,\sigma_2)}}$ be the union of the $K$-cubes in $
Y^{\sigma_1,\sigma_2}_{A^{(\sigma_1,\sigma_2)}}$ such that for each $
B\subset
Y^{\sigma_1,\sigma_2}_{A^{(\sigma_1,\sigma_2)},\nu^{(\sigma_1,\sigma_2)}}$,
the $R^{1/2}$-cube intersecting $B$ contains $\approx
\nu^{(\sigma_1,\sigma_2)}$ cubes from $Y^{\sigma_1,\sigma_2}_{A^{(\sigma_1,
\sigma_2)}}$.
\end{enumerate}

The dyadic numbers $A^{(\sigma _{1},\sigma _{2})}$ and $\nu ^{(\sigma
_{1},\sigma _{2})}$ making significant contribution can be assumed to be
between $R^{-C}$ and $R^{C}$. Therefore, there exist some dyadic numbers $
A^{(\sigma _{1},\sigma _{2})}$ and $\nu ^{(\sigma _{1},\sigma _{2})}$ such
that
\begin{equation*}
\sharp \{B:B\subset Y_{A^{(\sigma _{1},\sigma _{2})},\nu ^{(\sigma
_{1},\sigma _{2})}}^{\sigma _{1},\sigma _{2}}\}=:N^{(\sigma _{1},\sigma
_{2})}\gtrsim \frac{N^{(0)}}{(\log K)^{2}(\log R)^{2}}.
\end{equation*}
Fix that choice of $A^{(\sigma _{1},\sigma _{2})},\nu ^{(\sigma _{1},\sigma
_{2})}$ and denote $Y_{A^{(\sigma _{1},\sigma _{2})},\nu ^{(\sigma
_{1},\sigma _{2})}}^{\sigma _{1},\sigma _{2}}$ by $Y^{\sigma _{1},\sigma
_{2}}$ for convenience. Then, in the $\Omega _{\sigma _{1},\sigma _{2}}$
-case we have that
\begin{equation*}
\Vert e^{it\phi (D)}f_{\Omega _{0}}\Vert _{L^{6}(Y^{0})}\lesssim
_{\varepsilon }R^{\varepsilon }\Vert e^{it\phi (D)}f_{\Omega _{\sigma
_{1},\sigma _{2}}}\Vert _{L^{6}(Y^{\sigma _{1},\sigma _{2}})},
\end{equation*}
and
\begin{equation*}
\frac{N^{(0)}}{(\log K)^{2}(\log R)^{2}}\lesssim N^{(\sigma _{1},\sigma
_{2})}\leq N^{(0)},\quad \nu ^{(\sigma _{1},\sigma _{2})}\leq \nu
^{(0)},\quad \lambda ^{(\sigma _{1},\sigma _{2})}\leq \lambda ^{(0)},
\end{equation*}
where
\begin{equation*}
\lambda ^{(\sigma _{1},\sigma _{2})}:=\max_{B^{3}(x^{\prime },r)\subset
Q_{R},x^{\prime }\in \mathbb{R}^{3},r\geq K}\frac{\sharp \{B_{k}\subset
Y^{\sigma _{1},\sigma _{2}}:B_{k}\subset B^{3}(x^{\prime },r)\}}{r^{\alpha }}
.
\end{equation*}
We are going to estimate $\Vert e^{it\phi (D)}f_{\Omega _{\sigma _{1},\sigma
_{2}}}\Vert _{L^{6}(Y^{\sigma _{1},\sigma _{2}})}$. By the triangle
inequality, one has
\begin{equation}
\Vert e^{it\phi (D)}f_{\Omega _{\sigma _{1},\sigma _{2}}}\Vert
_{L^{6}(B)}\lesssim \sum_{\tau \subset \Omega _{\sigma _{1},\sigma
_{2}}}\Vert e^{it\phi (D)}f_{\tau }\Vert _{L^{6}(B)}.  \label{trivaldec}
\end{equation}
For each $\tau $, we deal with $e^{it\phi (D)}f_{\tau }$ by the rescaling
trick from Case (c) in Subsection 2.2 of \cite{LiMiaoZheng2020}. To do it,
we need to further decompose $f_{\tau }$ in physical space and perform
dyadic pigeonholing several times. First we divide the physical square
$[0,R]^{2}$ into $\frac{R}{K^{\frac{1}{2}}\sigma _{1}^{-\frac{m-2}{2}}}
\times \frac{R}{K^{\frac{1}{2}}\sigma _{2}^{-\frac{m-2}{2}}}$-rectangles
$D$. For each pair $(\tau,D)$, let $f_{\Box _{\tau ,D}}$
be the function formed by cutting off $f$ on
the rectangle $D$ (with a Schwartz tail) in physical space and the rectangle
$\tau $ in Fourier space. Note that $e^{it\phi (D)}f_{\Box _{\tau ,D}}$ is
essentially supported on an
$\frac{R}{K^{\frac{1}{2}}\sigma _{1}^{-\frac{m-2}{2}}}
\times \frac{R}{K^{\frac{1}{2}}\sigma _{2}^{-\frac{m-2}{2}}}
\times R$-box, which is denoted by $\Box _{\tau ,D}$.
The long axis of $\Box _{\tau ,D}$ is parallel to the normal direction of
the surface $\Sigma _{0}$ at the left bottom corner of $\tau $.
For a fixed $\tau $, the different boxes $\Box _{\tau ,D}$
tile $Q_{R}$. We have
\begin{equation*}
f_{\Omega _{\sigma _{1},\sigma _{2}}}=\sum_{\tau \subset {\Omega _{\sigma
_{1},\sigma _{2}}}}\sum_{D}f_{\Box _{\tau ,D}},
\end{equation*}
and write
\begin{equation*}
f_{\Omega _{\sigma _{1},\sigma _{2}}}=\sum_{\Box }f_{\Box }
\end{equation*}
for simplicity. For each $\tau $, a given $K$-cube $B$ lies in exactly one
box $\Box _{\tau ,D}$. Recall that $K=R^{\delta }$, where $\delta
=\varepsilon ^{100}$. Denote
\begin{equation*}
R_{1}:=\frac{R}{K}=R^{1-\delta },\quad K_{1}=R_{1}^{\delta }=R^{\delta
-\delta ^{2}}.
\end{equation*}
Tile $\Box $ by $\sigma _{1}^{\frac{m-2}{2}}K^{1/2}K_{1}\times \sigma _{2}^{
\frac{m-2}{2}}K^{1/2}K_{1}\times KK_{1}$-tube $S$, and also tile $\Box $ by $
\frac{R^{1/2}}{\sigma _{2}^{-\frac{m-2}{2}}}\times \frac{R^{1/2}}{\sigma
_{2}^{-\frac{m-2}{2}}}\times K^{\frac{1}{2}}R^{\frac{1}{2}}$-tubes $
S^{\prime }$ (all running parallel to the long axis of $\Box $). After
rescaling, the $\Box $ becomes an $R_{1}$-cube, the tubes $S^{\prime }$ and $
S$ become lattice $R_{1}^{\frac{1}{2}}$-cubes and $K_{1}$-cubes, respectively.

We regroup tubes $S$ and $S^{\prime }$ inside each $\Box$ as follows:

\begin{enumerate}
[1)]
\item Sort those tubes $S$ which intersect $Y^0$ according to the value $
\Vert e^{it \phi(D)}f_{\Box}\Vert_{L^6(S)}$ and the number of $K$-cubes
contained in it. For dyadic numbers $\eta,\beta_1$, we use $\mathbb{S}
_{\Box,\eta,\beta_1}$ to stand for the collection of tubes $S\subset \Box$
each of which containing $\sim \eta$ $K$-cubes and $\Vert e^{it
\phi(D)}f_{\Box}\Vert_{L^6(S)}\sim \beta_1$.

\item For fixed $\eta,\beta_1$, we sort the tubes $S^{\prime }\subset \Box$
according to the number of the tubes $S\in \mathbb{S}_{\Box,\eta,\beta_1}$
contained in it. For dyadic number $\nu_1$, let $\mathbb{S}
_{\Box,\eta,\beta_1,\nu_1}$ be the subcollection of $\mathbb{S}
_{\Box,\eta,\beta_1}$ such that for each $S\in \mathbb{S}_{\Box,\eta,
\beta_1,\nu_1}$, the tube $S^{\prime }$ containing $S$ contains $\sim \nu_1$
tubes from $\mathbb{S}_{\Box,\eta,\beta_1}$.

\item For fixed $\eta ,\beta _{1},\nu _{1}$, we sort the boxes $\Box $
according to the value $\Vert f\Vert _{2}$, the number $\sharp \mathbb{S}
_{\Box ,\eta ,\beta _{1},\nu _{1}}$ and the value $\lambda _{1}$ defined
below. For dyadic numbers $\beta _{2},N_{1},\lambda _{1}$, let $\mathbb{B}
_{\eta ,\beta _{1},\nu _{1},\beta _{2},M_{1},\lambda _{1}}$ denote the
collection of boxes $\Box $ satisfying that
\begin{equation*}
\Vert f\Vert _{2}\sim \beta _{2},\quad \sharp \mathbb{S}_{\Box ,\eta ,\beta
_{1},\nu _{1}}\sim N_{1}
\end{equation*}
and
\begin{equation}
\max_{T_{r}\subset \Box :r\geq K_{1}}\frac{\sharp \{S\in \mathbb{S}_{\Box
,\eta ,\beta _{1},\nu _{1}}:S\subset T_{r}\}}{r^{2}}\sim \lambda _{1},
\label{4gamma11}
\end{equation}
where $T_{r}$ are $\sigma _{2}^{\frac{m-2}{2}}K^{\frac{1}{2}}r\times \sigma
_{2}^{\frac{m-2}{2}}K^{\frac{1}{2}}r\times Kr$-tubes in $\Box $ running
parallel to the long axis of $\Box $.
\end{enumerate}

Let $Y_{\Box ,\eta ,\beta _{1},\nu _{1}}^{\sigma _{1},\sigma _{2}}$ denote $
\{S:S\subset \mathbb{S}_{\Box ,\eta ,\beta _{1},\nu _{1}}\}$ and $\chi
_{Y_{\Box ,\eta ,\beta _{1},\nu _{1}}^{\sigma _{1},\sigma _{2}}}$ be the
corresponding characteristic function. Without loss of generality, we can assume
$\Vert f\Vert _{2}=1$. Therefore, there are only $O(logR)$ significant
choices for each dyadic number. By pigeonholing, we can choose $\eta ,\beta
_{1},\nu _{1},\beta _{2},N_{1},\lambda _{1}$ such that
\begin{equation}
\begin{split}
& \quad \Vert e^{it\phi (D)}f_{\Omega _{\sigma _{1},\sigma _{2}}}\Vert
_{L^{6}(B)} \\
& \lesssim (logR)^{6}K^{1/2}\Big(\sum_{\Box \in \mathbb{B}_{\eta ,\beta
_{1},\nu _{1},\beta _{2},N_{1},\lambda _{1}},B\subset Y_{\Box ,\eta ,\beta
_{1},\nu _{1}}}\Vert e^{it\phi (D)}f_{\Box }\Vert _{L^{6}(\omega _{B})}^{2}
\Big)^{1/2}
\end{split}
\label{4decoup314}
\end{equation}
holds for a fraction $\gtrsim (logR)^{-6}$ of all $K$-cubes $B$, where we
have used the fact that
\begin{equation*}
\sharp \{\tau :\tau \subset \Omega _{\sigma _{1},\sigma _{2}}\}\lesssim K.
\end{equation*}
For brevity, we denote by
\begin{equation*}
Y_{\Box }^{\sigma _{1},\sigma _{2}}:=Y_{\Box ,\eta ,\beta _{1},\nu
_{1}}^{\sigma _{1},\sigma _{2}},\quad \mathbb{B}:=\mathbb{B}_{\eta ,\beta
_{1},\nu _{1},\beta _{2},N_{1},\lambda _{1}}.
\end{equation*}
Finally, we sort the $K$-cubes $B$ satisfying \eqref{4decoup314} by $\sharp
\{\Box \in \mathbb{B}:B\subset Y_{\Box }^{\sigma _{1},\sigma _{2}}\}$. Let $
Y^{\prime \sigma _{1},\sigma _{2}}$ be a union of $K$-cubes $B$ obeying
\begin{equation}
\Vert e^{it\phi (D)}f_{\Omega _{\sigma _{1},\sigma _{2}}}\Vert
_{L^{6}(B)}\lesssim (logR)^{6}K^{1/2}\Big(\sum_{\Box \in \mathbb{B},B\subset
Y_{\Box }^{0}}\Vert e^{it\phi (D)}f_{\Box }\Vert _{L^{6}(\omega _{B})}^{2}
\Big)^{1/2}  \label{4decoup315}
\end{equation}
and
\begin{equation}
\sharp \{\Box \in \mathbb{B}:B\subset Y_{\Box }^{\sigma _{1},\sigma
_{2}}\sim \mu \}  \label{4boxdyadic}
\end{equation}
for some dyadic number $1\leq \mu \leq K^{O(1)}$. Moreover, the number of $K$
-cubes $B$ in $Y^{\prime }$ is $\gtrsim (logR)^{-7}N$.

By our assumption that $\Vert e^{it
\phi(D)}f_{\Omega_{\sigma_1,\sigma_2}}\Vert_{L^6(B_k)}$ is essentially
constant in $k=1,2,\cdots,N^0$, we have
\begin{equation}  \label{4ineq6}
\Vert e^{it \phi(D)}f_{\Omega_{\sigma_1,\sigma_2}}
\Vert^6_{L^6(Y^{\sigma_1,\sigma_2})} \lesssim K^3 (logR)^7 \sum_{B\subset
Y^{\prime }} \Vert e^{it \phi(D)}f_{\Omega_{\sigma_1,\sigma_2}}
\Vert^6_{L^6(B)}.
\end{equation}
For each $B\subset Y^{\prime }$, it follows from \eqref{4decoup315},
\eqref{4boxdyadic} and H\"{o}lder's inequality that
\begin{equation}  \label{4decoup316}
\Vert e^{it \phi(D)}f_{\Omega_{\sigma_1,\sigma_2}} \Vert^6_{L^6(B)} \lesssim
(logR)^{36} K^3 \mu^2 \Big(\sum_{\Box \in \mathbb{B},B\subset
Y^{\sigma_1,\sigma_2}_{\Box}}\Vert e^{it \phi(D)}f_{\Box}
\Vert^6_{L^6(\omega_{B})} \Big)^{1/2}.
\end{equation}
Combining \eqref{4ineq6} with \eqref{4decoup316}, one has
\begin{equation}  \label{4decoup7}
\Vert e^{it \phi(D)}f_{\Omega_{\sigma_1,\sigma_2}}
\Vert_{L^6(Y^{\sigma_1,\sigma_2})} \lesssim (logR)^{C} K^{1/2}\mu^{1/3}\Big(
\sum_{\Box \in \mathbb{B}}\Vert e^{it \phi(D)}f_{\Box}
\Vert^6_{L^6(Y^{\sigma_1,\sigma_2}_{\Box})} \Big)^{1/6}.
\end{equation}
Next, we apply rescaling to each $\Vert e^{it \phi(D)}f_{\Box}
\Vert_{L^6(Y^{\sigma_1,\sigma_2}_{\Box})}.$ For each $\sigma_1^{-\frac{m-2}{2
}}K^{-1/2}\times \sigma_2^{-\frac{m-2}{2}}K^{-1/2}$-rectangle $\tau=
\tau_{\Box}$ in $\Omega_{\sigma_1,\sigma_2}$, we write
\begin{equation*}
\left\{
\begin{array}{l}
\xi_1 = a_1 + \sigma_1^{-\frac{m-2}{2}}K^{-1/2}\eta_1, \\
\xi_2 = a_2 +\sigma_2^{-\frac{m-2}{2}}K^{-1/2}\eta_2.
\end{array}
\right.
\end{equation*}
Then
\begin{equation*}
\vert e^{it \phi(D)}f_{\Box}(x)\vert= \sigma_1^{-\frac{m-2}{4}}\sigma_2^{-
\frac{m-2}{4}}K^{-1/2}\vert e^{i\tilde{t} \tilde{\phi}(D)}g(\tilde{x})\vert
\end{equation*}
with $\Vert g \Vert_2= \Vert f \Vert_2$ and
\begin{equation*}
\hat{g}(\eta_1,\eta_2):=\sigma_1^{-\frac{m-2}{4}}\sigma_2^{-\frac{m-2}{4}
}K^{-1/2}\hat{f}(a_1+\sigma_1^{-\frac{m-2}{2}}K^{-1/2}\eta_1, a_2+\sigma_2^{-
\frac{m-2}{2}}K^{-1/2}\eta_2),
\end{equation*}
where
\begin{equation*}
\begin{split}
& \quad \tilde{\phi}(\eta_1,\eta_2) \\
&:=\sum_{i=1}^2\Big(\frac{m(m-1)}{2}a_i^{m-2}\sigma_i^{-(m-2)}\eta^2_i \\
&\quad +\frac{m(m-1)(m-2)}{6}a_i^{m-3}\sigma_i^{-\frac{3(m-2)}{2}
}K^{-1/2}\eta^3_i+\dots +\sigma_i^{-\frac{m(m-2)}{2}}K^{-\frac{m-2}{2}
}\eta^m_i\Big),
\end{split}
\end{equation*}
and
\begin{equation*}
\left\{
\begin{array}{l}
\tilde{x}_1 := \sigma_1^{-\frac{m-2}{2}}K^{-1/2}(x_2+ma_1^{m-1}t), \\
\tilde{x}_2 := \sigma_2^{-\frac{m-2}{2}}K^{-1/2}(x_2+ma_2^{m-1}t), \\
\tilde{t} := K^{-1}t.
\end{array}
\right.
\end{equation*}
For brevity, we denote the above relation by $(\tilde{x},\tilde{t})=\mathcal{
L}_{0}(x,t)$. Therefore, we have
\begin{equation}  \label{4relation}
\Vert e^{it \phi(D)}f_{\Box}
\Vert_{L^6(Y^{\sigma_1,\sigma_2}_{\Box})}=\rho^{-\frac{m-2}{6}}K^{-1/6}\Vert
e^{i\tilde{t} \tilde{\phi}(D)}g(\tilde{x}) \Vert_{L^6(\tilde{Y}
^{\sigma_1,\sigma_2})},
\end{equation}
where $\tilde{Y}^{\sigma_1,\sigma_2}=\mathcal{L}_{0}(Y^{\sigma_1,\sigma_2}_{
\Box})$. Since $a_i \sim \sigma_i$ and $\frac{1}{2} \geq \sigma_i \geq
K^{-1/m}$, it is easy to check that $\tilde{\phi}$ is a phase function of
elliptic type, namely,
\begin{equation*}
\partial^2_{\eta_i}\tilde{\phi}\sim_m 1,\quad i=1,2,
\end{equation*}
and
\begin{equation*}
\vert \partial^l_{\eta_2}\tilde{\phi} \vert \lesssim_m 1,\quad 3\leq l \leq
m.
\end{equation*}
Applying Proposition 3.1 in \cite{DuZhang19-Annals} to the term on the
right-hand side of \eqref{4relation} at scale $R_1$, we deduce
\begin{equation}  \label{4induction}
\Vert e^{it \phi(D)}f_{\Box}
\Vert_{L^6(Y^{\sigma_1,\sigma_2}_{\Box})}\lesssim \sigma_1^{-\frac{m-2}{6}
}\sigma_2^{-\frac{m-2}{6}} K^{-1/6}N_1^{-\frac{1}{3}}\lambda_1^{\frac{1}{6}
}\nu_1^{\frac{1}{6}}(\frac{R}{K})^{\frac{\alpha}{12}+\varepsilon}\Vert
f_{\Box}\Vert_2.
\end{equation}
Using a similar argument as in the proof of inequality (3.24) in \cite
{DuZhang19-Annals}, one has
\begin{equation}  \label{4ineq8}
\frac{\mu}{\sharp \mathbb{B}} \lesssim \frac{(logR)^7 N_1 \eta}{N^0}
\end{equation}
and
\begin{equation*}
\lambda_1 \eta \lesssim \max_{T_r\subset \Box: r\geq K_1}\frac{\sharp
\{B\subset Y^0: B\subset T_r\}}{r^{\alpha}}\leq \frac{
\lambda^{(0)}(K^{O(1)}r)^{{\alpha}}}{r^{\alpha}},
\end{equation*}
where we have used the fact that one can cover a $\sigma_1^{\frac{m-2}{2}
}K^{1/2}r \times \sigma_2^{\frac{m-2}{2}}K^{1/2}r \times Kr$-tube $T_r$ by $
K^{O(1)}$ finitely overlapping $(\min\{\sigma_1,\sigma_2\})^{\frac{m-2}{2}
}K^{1/2}r$-balls. Hence, we get
\begin{equation}  \label{4ineq9}
\eta \lesssim \frac{\lambda^0 K^{O(1)}}{\lambda_1}.
\end{equation}

Now we relate $\nu _{1}$ and $\nu ^{(0)}$ by considering the number of $K$
-cubes in each relevant $\frac{R^{1/2}}{\sigma _{1}^{-\frac{m-2}{2}}}\times
\frac{R^{1/2}}{\sigma _{2}^{-\frac{m-2}{2}}}\times K^{\frac{1}{2}}R^{\frac{1
}{2}}$-tube $S^{\prime }$. Recall that each relevant $S^{\prime }$ contains $
\sim \nu _{1}$ tubes $S$ in $Y_{\Box }^{0}$ and each such $S$ contains $\sim
\eta $ $K$-cubes. On the other hand, we can cover $S^{\prime }$ by $K^{O(1)}$
finitely overlapping $R^{\frac{1}{2}}$-cubes, and each $R^{
\frac{1}{2}}$-cube contains $\lesssim \nu ^{(0)}$ many $K$-cubes in $Y^{0}$.
Thus, it follows that
\begin{equation}
\nu _{1}\lesssim \frac{K^{\frac{1}{2}}\nu ^{(0)}}{\eta }.  \label{4ineq10}
\end{equation}
By inserting \eqref{4ineq8}, \eqref{4ineq9} and \eqref{4ineq10} into
\eqref{4decoup7}, we derive
\begin{equation}
\begin{split}
& \quad \Vert e^{it\phi (D)}f_{\Omega _{\sigma _{1},\sigma _{2}}}\Vert
_{L^{6}(Y^{\sigma _{1},\sigma _{2}})} \\
& \leq C_{\varepsilon }K^{O(1)}\sigma _{1}^{-\frac{m-2}{6m}}\sigma _{2}^{-
\frac{m-2}{6m}}(N^{(0)})^{-\frac{1}{3}}(\nu ^{(0)})^{\frac{1}{6}}(\lambda
^{(0)})^{\frac{1}{6}}R^{\frac{\alpha }{12}+\varepsilon }\Vert f\Vert _{2}.
\end{split}
\label{4omega0}
\end{equation}
Recall that $\sigma _{i}\geq K^{-1/m}$. Combining \eqref{4omega0} with
\eqref{4decompose0}, we obtain
\begin{equation*}
\Vert e^{it\phi (D)}f_{\Omega _{0}}\Vert _{L^{6}(Y^{0})}\leq C_{\varepsilon
}K^{O(1)}(N^{(0)})^{-\frac{1}{3}}(\nu ^{(0)})^{\frac{1}{6}}(\lambda ^{(0)})^{
\frac{1}{6}}R^{\frac{\alpha }{12}+\varepsilon }\Vert f\Vert _{2},
\end{equation*}
which verifies inequality \eqref{omega0}.

For the $\Omega_3$-case, we sort the $K$-cubes in $Y^3$ as follows:

\begin{enumerate}
[(1)]
\item For any dyadic number $A^{(3)}$, let $Y^3_{A^{(3)}}$ be the union of
the $K$-cubes in $Y^3$ satisfying
\begin{equation}  \label{p6line40}
\Vert e^{it \phi(D)}f_{\Omega_3} \Vert_{L^6(B)}\approx A^{(3)}.
\end{equation}

\item Fix $A^{(3)}$, for any dyadic number $\nu^{(3)}$, let $
Y^3_{A^{(3)},\nu^{(3)}}$ be the union of the $K$-cubes in $Y^3_{A^{(3)}}$
such that for each $B\subset Y^3_{A^{(3)},\nu^{(3)}}$, the $R^{1/2}$-cube
intersecting $B$ contains $\approx \nu^{(3)}$ cubes from $Y^3_{A^{(3)}}$.
\end{enumerate}

The dyadic numbers $A^{(3)}$ and $\nu^{(3)}$ making significant contribution
can be assumed to be between $R^{-C}$ and $R^C$. Therefore, there exist some
dyadic numbers $A^{(3)}$ and $\nu^{(3)}$ such that
\begin{equation*}
\sharp \{B: B \subset Y^3_{A^{(3)},\nu^{(3)}}\} =:N^{(3)} \gtrsim \frac{N}{
(\log R)^2}.
\end{equation*}
Fix that choice of $A^{(3)}, \nu^{(3)}$ and denote $Y^3_{A^{(3)},\nu^{(3)}}$
by $Y^3$ for convenience. Then in the $\Omega_3$-case, we have
\begin{equation}  \label{p6line52}
\Vert e^{it \phi(D)}f \Vert_{L^6(Y)}\lesssim_{\varepsilon}R^{\varepsilon}
\Vert e^{it \phi(D)}f_{\Omega_3} \Vert_{L^6(Y^3)},
\end{equation}
and
\begin{equation*}
\frac{N}{(\log R)^2} \lesssim N^{(3)} \leq N,\quad \nu^{(3)} \leq \nu,\quad
\lambda^{(3)} \leq \lambda,
\end{equation*}
where
\begin{equation*}
\lambda^{(3)}:=\max_{B^3(x^{\prime },r)\subset Q_R,x^{\prime }\in \mathbb{R}
^3,r\geq 1}\frac{\sharp\{B_k \subset Y^3: B_k \subset B^3(x^{\prime },r)\}}{
r^{\alpha}}.
\end{equation*}
We are going to estimate $\Vert e^{it \phi(D)}f_{\Omega_3} \Vert_{L^6(Y^3)}$
. We deal with $e^{it\phi(D)}f_{\Omega_3}$ by rescaling and induction on
scales. To do it, we further decompose $f_{\Omega_3}$ in physical space and
perform dyadic pigeonholing several times.

Firstly, we divide the physical square $[0,R]^2$ into $\frac{R}{K^{1-\frac{1
}{m}}}\times \frac{R}{K^{1-\frac{1}{m}}}$-rectangles $D$. For each pair $D$,
let $f_{\Box_{\Omega_3,D}}$ be the function formed by cutting off $
f_{\Omega_3}$ on the rectangle $D$ (with a Schwartz tail) in physical space.
Thus, we have
\begin{equation*}
f_{\Omega_3}=\sum_{D}f_{\Box_{\Omega_3,D}},
\end{equation*}
and write $f_{\Omega_3}=\sum_{\Box}f_{\Box}$ for abbreviation. Note that
$e^{it \phi(D)}f_{\Box_{\Omega_3,D}}$ is essentially supported on
an $\frac{R}{K^{1-\frac{1}{m}}}\times \frac{R}{K^{1-\frac{1}{m}}}\times R$
-box, which is denoted by $\Box_{\Omega_3,D}$. The long axis of $
\Box_{\Omega_3,D}$ is parallel to the vector $(0,0,1)$. The different boxes
$\Box_{\Omega_3,D}$ tile $Q_R$. In particular, a given $K$-cube $B$ lies
in exactly one box $\Box_{\Omega_3,D}$.
For each $K$-cube $B$, it holds trivially that
\begin{equation}  \label{trivialL6B}
\Vert e^{it\phi(D)}f_{\Omega_3} \Vert_{L^6(B)}\lesssim \Big(\sum_{\Box}
\Vert e^{it\phi(D)}f_{\Box} \Vert^6_{L^6(B)} \Big)^{1/6},
\end{equation}
and
\begin{equation}  \label{case3box}
\sharp \{\Box\in \mathbb{B}: B\subset Y^3_{\Box}\}\sim 1.
\end{equation}
Recall that $K=R^{\delta}$, where $\delta=\varepsilon^{100}$. Denote
\begin{equation*}
R_1:=\frac{R}{K}=R^{1-\delta},\quad K_1=R^{\delta}_1=R^{\delta-\delta^2}.
\end{equation*}
Tile $\Box$ by $K^{\frac{1}{m}}K_1\times K^{\frac{1}{m}}K_1 \times K K_1$
-tube $S$, and also tile $\Box$ by $\frac{R^{1/2}}{K^{\frac{1}{2}-\frac{1}{m}
}}\times \frac{R^{1/2}}{K^{\frac{1}{2}-\frac{1}{m}}}\times K^{\frac{1}{2}}R^{
\frac{1}{2}}$-tubes $S^{\prime }$ (all running parallel to the long axis of $
\Box$). After rescaling the $\Box$ becomes an $R_1$-cube, the tubes $
S^{\prime }$ and $S$ becomes lattice $R_1^{\frac{1}{2}}$-cubes and $K_1$
-cubes, respectively. We regroup the tubes $S$ and $S^{\prime }$ inside each
$\Box$ as follows:

\begin{enumerate}
[(1)]
\item Sort those tubes $S$ which intersect $Y^{3}$ according to the value $
\Vert e^{it\phi (D)}f_{\Box }\Vert _{L^{6}(S)}$ and the number of $K$-cubes
contained in it. For dyadic numbers $\eta ,\beta _{1}$, we use $\mathbb{S}
_{\Box ,\eta ,\beta _{1}}$ to stand for the collection of tubes $S\subset
\Box $ each of which contains $\sim \eta $ $K$-cubes and $\Vert e^{it\phi
(D)}f_{\Box }\Vert _{L^{6}(S)}\sim \beta _{1}$.

\item For fixed $\eta,\beta_1$, we sort the tubes $S^{\prime }\subset \Box$
according to the number of the tubes $S\in \mathbb{S}_{\Box,\eta,\beta_1}$
contained in it. For each dyadic number $\nu_3$, let $\mathbb{S}
_{\Box,\eta,\beta_1,\nu_3}$ be the subcollection of $\mathbb{S}
_{\Box,\eta,\beta_1}$ such that for each $S\in \mathbb{S}_{\Box,\eta,
\beta_1,\nu_3}$, the tube $S^{\prime }$ containing $S$ contains $\sim \nu_3$
tubes from $\mathbb{S}_{\Box,\eta,\beta_1}$.

\item For fixed $\eta,\beta_1,\nu_3$, we sort the boxes $\Box$ according to
the value $\Vert f_{\Box} \Vert_2$, the number $\sharp \mathbb{S}
_{\Box,\eta,\beta_1,\nu_3}$ and the value $\lambda_3$ defined below. For
dyadic numbers $\beta_2,N_3,\lambda_3$, let $\mathbb{B}_{\eta,\beta_1,\nu_3,
\beta_2,N_3,\lambda_3}$ denote the collection of boxes $\Box$ each of which
satisfies
\begin{equation*}
\Vert f_{\Box} \Vert_2\sim \beta_2, \quad \sharp \mathbb{S}
_{\Box,\eta,\beta_1,\nu_3}\sim N_3
\end{equation*}
and
\begin{equation}  \label{gamma13}
\max_{T_r\subset \Box:r\geq K_1}\frac{\sharp\{S\in \mathbb{S}
_{\Box,\eta,\beta_1,\nu_3}:S\subset T_r\}}{r^{\alpha}}\sim\lambda_3,
\end{equation}
where $T_r$ are $K^{\frac{1}{m}}r\times K^{\frac{1}{m}}r \times Kr$-tubes in
$\Box$ running parallel to the long axis of $\Box$.
\end{enumerate}

Let $Y^3_{\Box,\eta,\beta_1,\nu_3}$ denote $\{S:S\subset \mathbb{S}
_{\Box,\eta,\beta_1,\nu_3}\}$. By \eqref{trivialL6B}, we have
\begin{equation*}
\Vert e^{it \phi(D)}f_{\Omega_3} \Vert^6_{L^6(B)}\lesssim
\sum_{\eta,\beta_1,\nu_3,\beta_2,N_3,\lambda_3}\Big(\sum_{\Box \in \mathbb{B}
_{\eta,\beta_1,\nu_3,\beta_2,N_3,\lambda_3},B\subset
Y^3_{\Box,\eta,\beta_1,\nu_3}}\Vert e^{it \phi(D)}f_{\Box}
\Vert^6_{L^6(\omega_{B})}\Big).
\end{equation*}
Thus, for each $B\subset Y^3$, we can choose $\eta,\beta_1,\nu_3,
\beta_2,N_3,\lambda_3$ depending on $B$ such that
\begin{equation*}
\Vert e^{it \phi(D)}f_{\Omega_3} \Vert^6_{L^6(B)} \lesssim (\log R)^6 \Big(
\sum_{\Box \in \mathbb{B}_{\eta,\beta_1,\nu_3,\beta_2,N_3,\lambda_3},B
\subset Y^3_{\Box,\eta,\beta_1,\nu_3}}\Vert e^{it \phi(D)}f_{\Box}
\Vert^6_{L^6(\omega_{B})}\Big).
\end{equation*}
Since there are only $O(\log R)$ significant choices for each dyadic number,
by pigeonholing, we can choose $\eta,\beta_1,\nu_3,\beta_2,N_3,\lambda_3$
such that
\begin{equation}  \label{case3314}
\begin{split}
& \quad \Vert e^{it \phi(D)}f_{\Omega_3} \Vert^6_{L^6(B)} \\
& \lesssim (\log R)^6 \Big(\sum_{\Box \in \mathbb{B}_{\eta,\beta_1,\nu_3,
\beta_2,N_3,\lambda_3},B\subset Y^3_{\Box,\eta,\beta_1,\nu_3}}\Vert e^{it
\phi(D)}f_{\Box} \Vert^6_{L^6(\omega_{B})}\Big)
\end{split}
\end{equation}
holds for a fraction $\gtrsim (\log R)^{-6}$ of all $K$-cubes $B$ in $Y^3$.
For brevity, we denote by
\begin{equation*}
Y^3_{\Box}:=Y^3_{\Box,\eta,\beta_1,\nu_3}, \quad \mathbb{B}:=\mathbb{B}
_{\eta,\beta_1,\nu_3,\beta_2,N_3,\lambda_3}.
\end{equation*}
Let $Y^{\prime}$ be the union of $K$-cubes $B$ each of which obeys
\begin{equation}  \label{case3315}
\Vert e^{it \phi(D)}f_{\Omega_3} \Vert_{L^6(B)} \lesssim \log R \Big(
\sum_{\Box \in \mathbb{B},B\subset Y^3_{\Box}}\Vert e^{it \phi(D)}f_{\Box}
\Vert^6_{L^6(B)}\Big)^{1/6}.
\end{equation}
The number of $K$-cubes $B$ in $Y^{\prime }$ is $\gtrsim (\log
R)^{-6}N^{(3)} $. This together with \eqref{p6line40} yields
\begin{equation}  \label{case3ineq6}
\Vert e^{it \phi(D)}f_{\Omega_3} \Vert^6_{L^6(Y^3)} \lesssim (\log
R)^{C}\sum_{B\subset Y^{\prime }} \Vert e^{it \phi(D)}f_{\Omega_3}
\Vert^6_{L^6(B)}.
\end{equation}
Putting \eqref{case3ineq6} and \eqref{case3315} together, one has
\begin{equation}  \label{case37}
\Vert e^{it \phi(D)}f_{\Omega_3} \Vert_{L^6(Y^3)} \lesssim (\log R)^{C}\Big(
\sum_{\Box \in \mathbb{B}}\Vert e^{it \phi(D)}f_{\Box}
\Vert^6_{L^6(Y^3_{\Box})} \Big)^{1/6}.
\end{equation}

Next, we apply rescaling to each $\Vert e^{it\phi (D)}f_{\Box }\Vert
_{L^{6}(Y_{\Box }^{3})}$ and induction on scales. Taking the change of
variables
\begin{equation*}
\left\{
\begin{array}{l}
\xi _{1}=K^{-\frac{1}{m}}\eta _{1}, \\
\xi _{2}=K^{-\frac{1}{m}}\eta _{2},
\end{array}
\right.
\end{equation*}
we have
\begin{equation*}
|e^{it\phi (D)}f_{\Box }(x)|=K^{-\frac{1}{m}}|e^{i\tilde{t}\phi (D)}g(\tilde{
x})|,
\end{equation*}
where
\begin{equation*}
\hat{g}(\eta _{1},\eta _{2}):=K^{-\frac{1}{m}}\hat{f_{\Box }}(K^{-\frac{1}{m}
}\eta _{1},K^{-\frac{1}{m}}\eta _{2}),\quad \Vert g\Vert _{2}=\Vert f_{\Box
}\Vert _{2},
\end{equation*}
and
\begin{equation*}
\left\{
\begin{array}{l}
\tilde{x}_{1}:=K^{-\frac{1}{m}}x_{1}, \\
\tilde{x}_{2}:=K^{-\frac{1}{m}}x_{2}, \\
\tilde{t}:=K^{-1}t.
\end{array}
\right.
\end{equation*}
For brevity, we denote the above relation by $(\tilde{x},\tilde{t})=\mathcal{
L}_{3}(x,t)$. Therefore, we have
\begin{equation}
\Vert e^{it\phi (D)}f_{\Box }\Vert _{L^{6}(Y_{\Box }^{3})}=K^{\frac{1}{6}-
\frac{2}{3m}}\Vert e^{i\tilde{t}\phi (D)}g(\tilde{x})\Vert _{L^{6}(\tilde{Y}
^{3})},  \label{case3relation}
\end{equation}
where $\tilde{Y}^{3}=\mathcal{L}_{3}(Y_{\Box }^{3})$. Hence, by
\eqref{case3relation} and the inductive hypothesis at scale $R_{1}$, we
deduce
\begin{equation}
\Vert e^{it\phi (D)}f_{\Box }\Vert _{L^{6}(Y_{\Box }^{3})}\lesssim K^{\frac{1
}{6}-\frac{2}{3m}}N_{3}^{-\frac{1}{3}}\lambda _{3}^{\frac{1}{6}}\nu _{3}^{
\frac{1}{6}}(\tfrac{R}{K})^{\frac{5}{12}-\frac{5-\alpha }{6m}+\varepsilon
}\Vert f_{\Box }\Vert _{2}.  \label{case3induction}
\end{equation}
By inserting \eqref{case3induction} into \eqref{case37}, we get
\begin{equation}
\Vert e^{it\phi (D)}f_{\Omega _{3}}\Vert_{L^{6}(Y^{3})}\lesssim
(N^{(3)})^{-\frac{1}{3}}(\nu^{(3)})^{\frac{1}{6}}
(\lambda^{(3)})^{\frac{1}{6}}
(\frac{R}{K})^{\frac{5}{12}-\frac{5-\alpha }{6m}+\varepsilon }
\big(\sum_{\Box \in \mathbb{B}}\Vert f_{\Box }\Vert_{2}^{6}\big)^{1/6}.
\label{case3inductionbox}
\end{equation}
Note that $\Vert f_{\Box }\Vert _{2}$ is essentially constant for $\Box \in
\mathbb{B}$. It follows that
\begin{equation}
\Vert e^{it\phi (D)}f_{\Omega _{3}}\Vert _{L^{6}(Y^{3})}\lesssim (\sharp
\mathbb{B})^{-1/3}(N^{(3)})^{-\frac{1}{3}}(\nu ^{(3)})^{\frac{1}{6}}(\lambda
^{(3)})^{\frac{1}{6}}R^{\frac{5}{12}-\frac{5-\alpha }{6m}+\varepsilon }\Vert
f\Vert _{2}.  \label{case3B}
\end{equation}
Consider the cardinality of the set $\{(\Box ,B):\Box \in \mathbb{B}
,B\subseteq Y_{\Box }^{3}\cap Y^{\prime }\}$. It is easy to see that there
is a lower bound
\begin{equation*}
\sharp \{(\Box ,B):\Box \in \mathbb{B},B\subseteq Y_{\Box }^{3}\cap
Y^{\prime -6}N^{(3)}.
\end{equation*}
On the other hand, by our choices of $N_{3}$ and $\eta $, for each $\Box \in
\mathbb{B}$, $Y_{\Box }^{3}$ contains $\sim N_{3}$ tubes $S$ and each $S$
contains $\sim \eta $ cubes in $Y^{3}$, so
\begin{equation*}
\sharp \{(\Box ,B):\Box \in \mathbb{B},B\subseteq Y_{\Box }^{3}\cap
Y^{\prime }\}\lesssim (\sharp \mathbb{B})N_{3}\eta .
\end{equation*}
Therefore, we get
\begin{equation}
\frac{1}{\sharp \mathbb{B}}\lesssim \frac{(\log R)^{6}N_{3}\eta }{N^{(3)}}.
\label{case3ineq8}
\end{equation}
Then by our choices of $\lambda _{3}$ as in \eqref{gamma13} and $\eta $, we
have
\begin{equation*}
\begin{split}
\lambda _{3}\eta & \sim \max_{T_{r}\subset \Box :r\geq K_{1}}\frac{\sharp
\{S:S\subseteq Y_{\Box }^{3}\cap T_{r}\}}{r^{\alpha }}\cdot \sharp
\{B:B\subseteq S\cap Y^{3}\} \\
& \lesssim \max_{T_{r}\subset \Box :r\geq K_{1}}\frac{\sharp \{B\subset
Y^{3}:B\subset T_{r}\}}{r^{\alpha }} \\
& \leq \frac{K^{1-\frac{1}{m}}\lambda ^{(3)}(K^{\frac{1}{m}}r)^{\alpha }}{
r^{\alpha }},
\end{split}
\end{equation*}
where the last inequality follows from the fact that we can cover a $K^{
\frac{1}{m}}r\times K^{\frac{1}{m}}r\times Kr$-tube $T_{r}$ by $\sim K^{1-
\frac{1}{m}}$ finitely overlapping $K^{\frac{1}{m}}r$-balls. Hence, we
obtain
\begin{equation}
\eta \lesssim \frac{\lambda ^{(3)}K^{1+\frac{\alpha -1}{m}}}{\lambda _{3}}.
\label{case3ineq9}
\end{equation}

Finally, we relate $\nu _{3}$ and $\nu ^{(3)}$ by considering the number of $
K$-cubes in each relevant $\frac{R^{1/2}}{K^{\frac{1}{2}-\frac{1}{m}}}\times
\frac{R^{1/2}}{K^{\frac{1}{2}-\frac{1}{m}}}\times K^{\frac{1}{2}}R^{\frac{1}{
2}}$-tube $S^{\prime }$. Recall that each relevant $S^{\prime }$ contains $
\sim \nu _{3}$ tubes $S$ in $Y_{\Box }$ and each such $S$ contains $\sim
\eta $ $K$-cubes. On the other hand, we can cover $S^{\prime }$ by $\sim
K^{1/2}$ finitely overlapping $R^{\frac{1}{2}}$-cubes, and each
$R^{\frac{1}{2}}$-cube contains $\lesssim \nu ^{(3)}$ many $K$-cubes in $
Y^{3}$. Thus, it follows that
\begin{equation}
\nu _{3}\lesssim \frac{K^{\frac{1}{2}}\nu ^{(3)}}{\eta }.
\label{case3ineq10}
\end{equation}
By inserting \eqref{case3ineq8}, \eqref{case3ineq9} and \eqref{case3ineq10}
into \eqref{case3B}, we derive
\begin{equation}
\Vert e^{it\phi (D)}f_{\Omega _{3}}\Vert _{L^{6}(Y^{3})}\lesssim
K^{-\varepsilon }(N^{(3)})^{-\frac{1}{3}}(\nu ^{(3)})^{\frac{1}{6}}(\lambda
^{(3)})^{\frac{1}{6}}R^{\frac{5}{12}-\frac{5-\alpha }{6m}+\varepsilon }\Vert
f\Vert _{2}.  \label{omega3}
\end{equation}
Since $K=R^{\delta }$ and $R$ can be assumed to be sufficiently large
compared to any constant depending on $\varepsilon $, we have
$K^{-\varepsilon }\ll 1$, and the induction closes. Recall that $\frac{N}{\log
R}\lesssim N^{(3)}\leq N,\nu ^{(3)}\leq \nu $ and $\lambda ^{(3)}\leq
\lambda $. It follows that
\begin{equation*}
\Vert e^{it\phi (D)}f_{\Omega _{3}}\Vert _{L^{6}(Y^{3})}\lesssim
_{\varepsilon }N^{-\frac{1}{3}}\lambda ^{\frac{1}{6}}\nu ^{\frac{1}{6}}R^{
\frac{5}{12}-\frac{5-\alpha }{6m}+\varepsilon }\Vert f\Vert _{2}.
\end{equation*}

For the $\Omega_1$-case, we will employ the following proposition, whose
proof is given in Section 4. Let $\psi(\xi_1,\xi_2):=\phi_1(\xi_1)+\xi^m_2$
be a class of smooth phase functions satisfying
\begin{equation}  \label{condition}
\phi_1^{(\prime\prime)}\sim 1;\quad\vert \phi_1^{(k)}\vert \lesssim 1,3\leq
k \leq m;\quad \phi_1^{(l)}=0,l\geq m+1
\end{equation}
on $[0,1]$.

\begin{proposition}
\label{mainlemma} For any $0< \varepsilon < \frac{1}{100}$, there exist
constants $C_{\varepsilon}$ and $\delta=\varepsilon^{100}$ such that the
following holds for all $R\geq 1$ and all $f$ with $supp\hat{f} \subset
[0,1]^2$. Suppose that $Y= \bigcup_{k=1}^{N}B_k$ is a union of lattice $K$
-cubes in $Q_R$ and each $R^{\frac{1}{2}}$-cube intersecting $Y$ contains $
\sim \nu$ many $K$-cubes in $Y$, where $K=R^{\delta}$. Suppose that $\Vert
e^{it \psi(D)}f \Vert_{L^6(B_k)}$ is essentially constant in $k=1,2,...,N$.
Given $\lambda$ by
\begin{equation*}  \label{gamm4}
\lambda:=\max_{B^3(x^{\prime },r)\subset Q_R,x^{\prime }\in \mathbb{R}
^3,r\geq K}\frac{\sharp\{B_k: B_k \subset B^3(x^{\prime },r)\}}{r^{\alpha}},
\end{equation*}
then
\begin{equation*}
\Vert e^{it \psi(D)}f \Vert_{L^6(Y)}\leq C_{\varepsilon}N^{-\frac{1}{3}
}\lambda^{\frac{1}{6}}\nu^{\frac{1}{6}}R^{\frac{1}{3}-\frac{4-\alpha}{6m}
+\varepsilon}\Vert f \Vert_2.
\end{equation*}
\end{proposition}

Assume that Proposition \ref{mainlemma} holds for a while, we estimate $
\Vert e^{it\phi(D)}f_{\Omega_1} \Vert_{L^6(Y^1)}$ as follows. We divide the
region $\Omega_1$ into a family of subregions
\begin{equation*}
\Omega_1=\bigcup_{\rho}\Omega_{1,\rho},
\end{equation*}
where $\Omega_{1,\rho}:= I_{\rho}\times [0,K^{-\frac{1}{m}}]$ with $
I_{\rho}=[\rho,2\rho]$. Here $\rho$ is a dyadic number in $[K^{-\frac{1}{m}},
\frac{1}{2}]$. We abbreviate $\Omega_{1,\rho}$ by $\Omega_{\rho}$ and write
\begin{equation*}
e^{it\phi(D)}f_{\Omega_1}=\sum_{\rho}e^{it\phi(D)}f_{\Omega_{\rho}}.
\end{equation*}
We will estimate $\Vert e^{it \phi(D)}f \Vert_{L^6(Y^1)}$ in several cases.
Loosely speaking, if the contribution from the region $\Omega_{\rho}$
dominates, we call it $\Omega_{\rho}$ case. More precisely, given a $K$-cube
$B$, by the triangle inequality, we have
\begin{equation*}
\Vert e^{it \phi(D)}f_{\Omega_{\rho}}\Vert_{L^6(B)}\gtrsim \frac{\Vert e^{it
\phi(D)}f_{\Omega_1}\Vert_{L^6(B)}}{\log K}
\end{equation*}
for at least one of the $\rho$'s. We sort the $K$-cubes in $Y^1$ as follows:
Denote
\begin{equation*}
\{B \subset Y^1: \Vert e^{it \phi(D)}f_{\Omega_{\rho}} \Vert_{L^6(B)}
\gtrsim \frac{\Vert e^{it \phi(D)}f_{\Omega_1}\Vert_{L^6(B)}}{\log K}\}
\end{equation*}
by $Y^{\rho}$. Clearly, one has $Y=\bigcup_{\rho}Y^{\rho}$. If $\sharp \{B:
B \subset Y^{\rho}\}\gtrsim \frac{N}{\log K}$, we call it $\Omega_{\rho}$
-case.

For the $\Omega_{\rho}$-case, we sort the $K$-cubes in $Y^{\rho}$ as follows:

\begin{enumerate}
[(1)]
\item For any dyadic number $A^{(\rho)}$, let $Y^{\rho}_{A^{(\rho)}}$ be the
union of the $K$-cubes in $Y^{\rho}$ satisfying
\begin{equation*}
\Vert e^{it \phi(D)}f_{\Omega_{\rho}} \Vert_{L^6(B)}\approx A^{(\rho)}.
\end{equation*}

\item Fix $A^{(\rho)}$, for any dyadic numbers $\nu^{(\rho)}$, let $
Y^{\rho}_{A^{(\rho)},\nu^{(\rho)}}$ be the union of the $K$-cubes in $
Y^{\rho}_{A^{(\rho)}}$ such that for each $B\subset
Y^{\rho}_{A^{(\rho)},\nu^{(\rho)}}$, the $R^{1/2}$-cube intersecting $B$
contains $\approx \nu^{(\rho)}$ cubes from $Y^{\rho}_{A^{(\rho)}}$.
\end{enumerate}

The dyadic numbers $A^{(\rho )},\nu ^{(\rho )}$ making significant
contributions can be assumed to be between $R^{-C}$ and $R^{C}$. Therefore,
there exist some dyadic numbers $A^{(\rho )},\nu ^{(\rho )}$ such that $
\sharp \{B:B\subset Y_{A^{(\rho )},\nu ^{(\rho )}}^{\rho }\}=:N^{(\rho
)}\gtrsim \frac{N}{\log K(\log R)^{4}}$ cubes $B$. Fix a choice of $A^{(\rho
)},\nu ^{(\rho )}$ and denote $Y_{A^{(\rho )},\nu ^{(\rho )}}^{\rho }$ by $
Y^{\rho }$ for convenience. Then, with a similar procedure as in the proof
of \eqref{c181}, in the $\Omega _{\rho }$-case we have
\begin{equation}
\Vert e^{it\phi (D)}f\Vert _{L^{6}(Y)}\lesssim _{\varepsilon }R^{\varepsilon
}\Vert e^{it\phi (D)}f_{\Omega _{\rho }}\Vert _{L^{6}(Y^{\rho })}.
\label{case1YandYrho}
\end{equation}
We also have
\begin{equation*}
\frac{N}{\log K(\log R)^{4}}\lesssim N^{(\rho )}\leq N,\quad \nu ^{(\rho
)}\leq \nu ,\quad \lambda ^{(\rho )}\leq \lambda ,
\end{equation*}
where
\begin{equation*}
\lambda ^{(\rho )}:=\max_{B^{3}(x^{\prime },r)\subset Q_{R},x^{\prime }\in
\mathbb{R}^{3},r\geq K}\frac{\sharp \{B_{k}\subset Y^{\rho }:B_{k}\subset
B^{3}(x^{\prime },r)\}}{r^{\alpha }}.
\end{equation*}
We are going to estimate $\Vert e^{it\phi (D)}f_{\Omega _{\rho }}\Vert
_{L^{6}(Y^{\rho })}$. For this purpose, we divide each $I_{\rho }$ further
into
\begin{equation*}
I_{\rho }=\bigcup_{j=1}^{\rho ^{\frac{m}{2}}K^{1/2}}I_{\rho ,j}.
\end{equation*}
Each $I_{\rho ,j}$ has length of $\rho ^{-\frac{m-2}{2}}K^{-1/2}$. We write
\begin{equation*}
\Omega _{\rho }=\bigcup \tau ,
\end{equation*}
where each $\tau $ in the above union has the form
\begin{equation*}
I_{\rho ,j}\times \lbrack 0,K^{-\frac{1}{m}}].
\end{equation*}
Therefore, we have
\begin{equation*}
e^{it\phi (D)}f_{\Omega _{\rho }}=\sum_{\tau \subset \Omega _{\rho
}}e^{it\phi (D)}f_{\tau }.
\end{equation*}
For each cube $B\subset Y^{\rho }$, we employ the following decoupling
inequality from \cite{Yang21}, which has root in the original paper of
Bourgain and Demeter \cite{BD15}.

\begin{lemma}
\label{decouplingrho} Suppose that $B$ is a $K$-cube contained in $Y^{\rho}$
. Then for any $\varepsilon >0$, there exists a constant $C_{\varepsilon}$
such that
\begin{equation*}  \label{decoupineqrho}
\Vert e^{it\phi(D)}f_{\Omega_{\rho}} \Vert_{L^6(B)}\leq
C_{\varepsilon}K^{\varepsilon}\Big(\sum_{\tau} \Vert e^{it\phi(D)}f_{\tau}
\Vert^2_{L^6(\omega_{B})} \Big)^{1/2},
\end{equation*}
where $\omega_{B}$ is a weight function essentially supported on $B$.
\end{lemma}

\begin{remark}
\label{omega1decouplingfornarrowball} The authors of \cite{DuZhang19-Annals}
use a broad-narrow analysis to deal with the Schr\"{o}dinger case. In
particular, they apply Bourgain-Demeter's decoupling inequality in \cite{BD15}
on each narrow ball. In the current paper, to estimate $\Vert
e^{it\phi(D)}f\Vert_{L^6(Y)}$, we sort the cubes in $Y$ into $
Y=\bigcup^3_{j=0}Y^{j}$, which can be regarded as a
``non-degenerate-degenerate" analysis. For the convenience of readers, in
this remark we call the cubes in $Y^0$ and $Y^1$ non-degenerate cubes and
degenerate cubes, respectively. Our strategy is to apply decoupling on each
degenerate cube, which has a similar flavor as the narrow ball in \cite
{DuZhang19-Annals}.
\end{remark}

For each $\tau$, we will handle $e^{it\phi(D)}f_{\tau}$ by rescaling and
induction on scales. To do it, we further decompose $f_{\Omega_{\rho}}$ in
physical space and perform dyadic pigeonholing several times.

Firstly, we divide the physical square $[0,R]^{2}$ into $\frac{R}{\rho ^{-
\frac{m-2}{2}}K^{1/2}}\times \frac{R}{K^{1-\frac{1}{m}}}$-rectangles $D$.
For each pair $(\tau ,D)$, let $f_{\Box _{\tau ,D}}$ be the function formed
by cutting off $f$ on the rectangle $D$ (with a Schwartz tail) in physical
space and the rectangle $\tau $ in Fourier space. We have
\begin{equation*}
f_{\Omega _{\rho }}=\sum_{\tau }\sum_{D}f_{\Box _{\tau ,D}},
\end{equation*}
and write $f_{\Omega _{\rho }}=\sum_{\Box }f_{\Box }$ for abbreviation.
Note that $e^{it\phi(D)}f_{\Box _{\tau ,D}}$ is essentially supported on an
$\frac{R}{\rho ^{-\frac{m-2}{2}}K^{1/2}}
\times \frac{R}{K^{1-\frac{1}{m}}}\times R$-box, which
is denoted by $\Box _{\tau ,D}$. The box $\Box _{\tau ,D}$ is parallel to
the normal direction of the surface $\Sigma _{\tau }$ at the left bottom
corner of $\tau $. For a fixed $\tau $, the different boxes $\Box _{\tau ,D}$
tile $Q_{R}$.  Note that a given $K$-cube $B$ lies in exactly one box
$\Box _{\tau ,D}$ for each $\tau $. By Lemma \ref{decouplingrho},
for each $K$-cube $B$, there holds
\begin{equation*}
\Vert e^{it\phi (D)}f_{\Omega _{\rho }}\Vert _{L^{6}(B)}\leq C_{\varepsilon
}K^{\varepsilon ^{2}}\Big(\sum_{\Box }\Vert e^{it\phi (D)}f_{\Box }\Vert
_{L^{6}(\omega _{B})}^{2}\Big)^{1/2}.
\end{equation*}
We denote
\begin{equation*}
R_{1}:=\frac{R}{K}=R^{1-\delta },\quad K_{1}=R_{1}^{\delta }=R^{\delta
-\delta ^{2}}.
\end{equation*}
Tile $\Box $ by $\rho ^{\frac{m-2}{2}}K^{\frac{1}{2}}K_{1}\times K^{\frac{1}{
m}}K_{1}\times KK_{1}$-tubes $S$, and also tile $\Box $ by $\frac{R^{1/2}}{
\rho ^{-\frac{m-2}{2}}}\times \frac{R^{1/2}}{K^{\frac{1}{2}-\frac{1}{m}}}
\times K^{\frac{1}{2}}R^{\frac{1}{2}}$-tubes $S^{\prime }$ (all running
parallel to the long axis of $\Box $). After rescaling the $\Box $ becomes an $
R_{1}$-cube, the tubes $S^{\prime }$ and $S$ become lattice $R_{1}^{\frac{1}{
2}}$-cubes and $K_{1}$-cubes, respectively. We regroup the tubes $S$ and $
S^{\prime }$ inside each $\Box $ as follows:

\begin{enumerate}
[(1)]
\item Sort those tubes $S$ which intersect $Y^{\rho}$ according to the value
$\Vert e^{it \phi(D)}f_{\Box}\Vert_{L^6(S)}$ and the number of $K$-cubes
contained in it. For dyadic numbers $\eta,\beta_1$, we use $\mathbb{S}
_{\Box,\eta,\beta_1}$ to stand for the collection of tubes $S\subset \Box$
each of which containing $\sim \eta$ $K$-cubes and $\Vert e^{it
\phi(D)}f_{\Box}\Vert_{L^6(S)}\sim \beta_1$.

\item For fixed $\eta,\beta_1$, we sort the tubes $S^{\prime }\subset \Box$
according to the number of the tubes $S\in \mathbb{S}_{\Box,\eta,\beta_1}$
contained in it. For dyadic number $\nu_{\rho}$, let $\mathbb{S}
_{\Box,\eta,\beta_1,\nu_{\rho}}$ be the sub-collection of $\mathbb{S}
_{\Box,\eta,\beta_1}$ such that for each $S\in \mathbb{S}_{\Box,\eta,
\beta_1,\nu_{\rho}}$, the tube $S^{\prime }$ containing $S$ contains $\sim
\nu_{\rho}$ tubes from $\mathbb{S}_{\Box,\eta,\beta_1}$.

\item For fixed $\eta ,\beta _{1},\nu _{\rho }$, we sort the boxes $\Box $
according to the value $\Vert f_{\Box }\Vert _{2}$, the number $\sharp
\mathbb{S}_{\Box ,\eta ,\beta _{1},\nu _{\rho }}$ and the value $\lambda
_{\rho }$ defined below. For dyadic numbers $\beta _{2},N_{\rho },\lambda
_{\rho }$, let $\mathbb{B}_{\eta ,\beta _{1},\nu _{\rho },\beta _{2},N_{\rho
},\lambda _{\rho }}$ denote the collection of boxes $\Box $ satisfying
\begin{equation*}
\Vert f_{\Box }\Vert _{2}\sim \beta _{2},\quad \sharp \mathbb{S}_{\Box ,\eta
,\beta _{1},\nu _{\rho }}\sim N_{\rho }
\end{equation*}
and
\begin{equation}
\max_{T_{r}\subset \Box :r\geq K_{1}}\frac{\sharp \{S\in \mathbb{S}_{\Box
,\eta ,\beta _{1},\nu _{\rho }}:S\subset T_{r}\}}{r^{\alpha }}\sim \lambda
_{\rho },  \label{gamma11rho}
\end{equation}
where $T_{r}$ are $K^{\frac{1}{2}}\rho ^{\frac{m-2}{2}}r\times K^{\frac{1}{m}
}r\times Kr$-tubes in $\Box $ running parallel to the long axis of $\Box $.
\end{enumerate}

Let $Y_{\Box ,\eta ,\beta _{1},\nu _{\rho}}^{\rho}$ denote $\{S:S\subset
\mathbb{S}_{\Box ,\eta ,\beta _{1},\nu _{\rho }}\}$. Thus, there are only
$O(\log R)$ significant choices for each dyadic number. By pigeonholing,
we can choose
$\eta ,\beta_{1},\nu _{\rho },\beta _{2},N_{\rho },\lambda _{\rho }$ so that
\begin{equation}
\begin{split}
& \quad \Vert e^{it\phi (D)}f_{\Omega _{\rho }}\Vert _{L^{6}(B)} \\
& \lesssim (\log R)^{6}K^{\varepsilon ^{2}}\Big(\sum_{\Box \in \mathbb{B}
_{\eta ,\beta _{1},\nu _{\rho },\beta _{2},N_{\rho },\lambda _{\rho
}},B\subset Y_{\Box ,\eta ,\beta _{1},\nu _{\rho }}}\Vert e^{it\phi
(D)}f_{\Box }\Vert _{L^{6}(\omega _{B})}^{2}\Big)^{1/2}
\end{split}
\label{decoup314rho}
\end{equation}
holds for a fraction $\gtrsim (\log R)^{-6}$ of all $K$-cubes $B$. For
brevity, we denote by
\begin{equation*}
Y_{\Box }^{\rho }:=Y_{\Box ,\eta ,\beta _{1},\nu _{\rho }}^{\rho },\quad
\mathbb{B}:=\mathbb{B}_{\eta ,\beta _{1},\nu _{\rho },\beta _{2},N_{\rho
},\lambda _{\rho }}.
\end{equation*}
We sort the $K$-cubes $B$ satisfying \eqref{decoup314rho} by $\sharp \{\Box
\in \mathbb{B}:B\subset Y_{\Box }^{\rho }\}$. Let $Y_{\mu }^{\prime }\subset
Y^{\rho }$ be a union of $K$-cubes $B$ that obey
\begin{equation}
\Vert e^{it\phi (D)}f_{\Omega _{\rho }}\Vert _{L^{6}(B)}\lesssim (\log
R)^{6}K^{\varepsilon ^{2}}\Big(\sum_{\Box \in \mathbb{B},B\subset Y_{\Box
}^{\rho }}\Vert e^{it\phi (D)}f_{\Box }\Vert _{L^{6}(\omega _{B})}^{2}\Big)
^{1/2}  \label{decoup315rho}
\end{equation}
and
\begin{equation}
\sharp \{\Box \in \mathbb{B}:B\subset Y_{\Box }^{\rho }\}\sim \mu .
\label{boxdyadicrho}
\end{equation}
Since the dyadic number $\mu $ must satisfy $1\leq \mu \leq R^{C}$, there
are $\sim \log R$ different choices for $\mu $, so by pigeonholing there
exists one such $\mu $ such that the number of cubes in $Y_{\mu }^{\prime }$
is $\gtrsim \frac{\sharp \{B:B\subset Y^{\rho }\}}{\log R}\gtrsim \frac{
N^{(\rho )}}{(\log R)^{7}}$. And then for simplicity we rename $Y_{\mu
}^{\prime }$ as $Y^{\prime }$. We have
\begin{equation}
\Vert e^{it\phi (D)}f_{\Omega _{\rho }}\Vert _{L^{6}(Y^{\rho })}^{6}\lesssim
(\log R)^{C}\sum_{B\subset Y^{\prime }}\Vert e^{it\phi (D)}f_{\Omega _{\rho
}}\Vert _{L^{6}(B)}^{6}.  \label{ineq6rho}
\end{equation}
For each $B\subset Y^{\prime }$, it follows from \eqref{decoup315rho},
\eqref{boxdyadicrho} and H\"{o}lder's inequality that
\begin{equation}
\Vert e^{it\phi (D)}f_{\Omega _{\rho }}\Vert _{L^{6}(B)}^{6}\lesssim (\log
R)^{C}K^{6\varepsilon ^{2}}\mu ^{2}\sum_{\Box \in \mathbb{B},B\subset
Y_{\Box }^{\rho }}\Vert e^{it\phi (D)}f_{\Box }\Vert _{L^{6}(\omega
_{B})}^{6}.  \label{decoup316rho}
\end{equation}
Putting \eqref{ineq6rho} and \eqref{decoup316rho} together, one has
\begin{equation}
\Vert e^{it\phi (D)}f_{\Omega _{\rho }}\Vert _{L^{6}(Y^{\rho })}\lesssim
(\log R)^{C}K^{\varepsilon ^{2}}\mu ^{1/3}\Big(\sum_{\Box \in \mathbb{B}
}\Vert e^{it\phi (D)}f_{\Box }\Vert _{L^{6}(Y_{\Box }^{\rho })}^{6}\Big)
^{1/6}.  \label{decoup7rho}
\end{equation}

Next, we apply rescaling to each $\Vert e^{it \phi(D)}f_{\Box}
\Vert_{L^6(Y^{\rho}_{\Box})}$ and run induction on scales. For each $\rho^{-
\frac{m-2}{2}}K^{-1/2}\times K^{-1/m}$-rectangle $\tau= \tau_{\Box}$ in $
\Omega_{\rho}$, we write
\begin{equation*}
\left\{
\begin{array}{l}
\xi_1 = a + \rho^{-\frac{m-2}{2}}K^{-\frac{1}{2}}\eta_1, \\
\xi_2 = K^{-\frac{1}{m}}\eta_2,
\end{array}
\right.
\end{equation*}
for some $\rho \leq a \leq 2\rho$. Then
\begin{equation*}
\vert e^{it \phi(D)}f_{\Box}(x)\vert= \rho^{-\frac{m-2}{4}}K^{-\frac{1}{4}-
\frac{1}{2m}}\vert e^{i\tilde{t}\psi(D)}g(\tilde{x})\vert,
\end{equation*}
with $\Vert g \Vert_2= \Vert f_{\Box} \Vert_2,$ and
\begin{equation*}
\hat{g}(\eta_1,\eta_2):= \rho^{-\frac{m-2}{4}}K^{-\frac{m+2}{4m}}\hat{
f_{\Box}}(a+K^{-\frac{1}{2}}\eta_1,K^{-\frac{1}{m}}\eta_2),
\end{equation*}
where
\begin{equation*}
\begin{split}
\psi(\eta_1,\eta_2) & :=\Big(\frac{m(m-1)}{2}a^{m-2}\rho^{-(m-2)}\eta^2_1 \\
&\quad +\frac{m(m-1)(m-2)}{6}a^{m-3}\rho^{-\frac{3(m-2)}{2}}K^{-1/2}\eta^3_1
\\
&\quad +...+\rho^{-\frac{m(m-2)}{2}}K^{-\frac{m-2}{2}}\eta^m_1\Big)+\eta^m_2
\\
&=:\phi_1(\eta_1)+\eta^m_2,
\end{split}
\end{equation*}
and
\begin{equation*}
\left\{
\begin{array}{l}
\tilde{x}_1 := \rho^{-\frac{m-2}{2}}K^{-\frac{1}{2}}(x_1+ma^{m-1}t), \\
\tilde{x}_2 := K^{-\frac{1}{m}}x_2, \\
\tilde{t} := K^{-1}t.
\end{array}
\right.
\end{equation*}
For brevity, we denote the above relation by $(\tilde{x},\tilde{t})=\mathcal{
L}_{1}(x,t)$. We write
\begin{equation}  \label{relationrho}
\Vert e^{it \phi(D)}f_{\Box} \Vert_{L^6(Y^{\rho}_{\Box})}=\rho^{-\frac{m-2}{6
}}K^{-\frac{1}{3m}}\Vert e^{i\tilde{t}\psi(D)}g(\tilde{x}) \Vert_{L^6(\tilde{
Y}^{\rho})},
\end{equation}
where $\tilde{Y}^{\rho}=\mathcal{L}_{1}(Y^{\rho}_{\Box})$. Note that $a\sim
\rho$ and $K^{-1/m}\leq \rho \leq \frac{1}{2}$. It is easy to see that the
new phase function $\psi(\eta_1,\eta_2)=\phi_1(\eta_1)+\eta^m_2$ satisfies
\begin{equation*}
\phi_1^{(\prime\prime)}\sim 1;\quad\vert \phi_1^{(k)}\vert \lesssim 1,3\leq
k \leq m;\quad \phi_1^{(l)}=0,l\geq m+1
\end{equation*}
on $[0,1]$. Therefore, applying Proposition \ref{mainlemma} to the term
\begin{equation*}
\Vert e^{i\tilde{t}\psi(D)}g(\tilde{x}) \Vert_{L^6(\tilde{Y}^{\rho})}
\end{equation*}
at scale $R_1$, we deduce
\begin{equation*}
\Vert e^{it \phi(D)}f_{\Box} \Vert_{L^6(Y^{\rho}_{\Box})}\lesssim \rho^{-
\frac{m-2}{6}}K^{-\frac{1}{3m}}N_{\rho}^{-\frac{1}{3}} \lambda_{\rho}^{\frac{
1}{6}}\nu_{\rho}^{\frac{1}{6}} (\tfrac{R}{K})^{\frac{1}{3}-\frac{4-\alpha}{6m
}+\varepsilon}\Vert f_{\Box} \Vert_2.
\end{equation*}
Note that $\frac{1}{3}-\frac{4-\alpha}{6m}\leq\frac{5}{12}-\frac{5-\alpha}{6m
}$ when $m \geq 2$. Consider the cardinality of the set $\{(\Box,B): \Box
\in \mathbb{B}, B\subseteq Y^{\rho}_{\Box}\cap Y^{\prime }\}$. By our choice
of $\mu$ as in \eqref{boxdyadicrho}, there is a lower bound
\begin{equation*}
\sharp \{(\Box,B): \Box \in \mathbb{B}, B\subseteq Y^{\rho}_{\Box}\cap
Y^{\prime -7}N^{({\rho})}\mu.
\end{equation*}
On the other hand, by our choices of $N_{\rho}$ and $\eta$, for each $\Box
\in \mathbb{B}$, $Y^{\rho}_{\Box}$ contains $\sim N_{\rho}$ tubes $S$ and
each $S$ contains $\sim \eta$ cubes in $Y^{\rho}$, so
\begin{equation*}
\sharp \{(\Box,B): \Box \in \mathbb{B}, B\subseteq Y^{\rho}_{\Box}\cap
Y^{\prime }\}\lesssim (\sharp \mathbb{B})N_{\rho}\eta.
\end{equation*}
Therefore, we get
\begin{equation}  \label{ineq8rho}
\frac{\mu}{\sharp \mathbb{B}} \lesssim \frac{(\log R)^7 N_{\rho} \eta}{N^{({
\rho})}}.
\end{equation}
Then by our choices of $\lambda_{\rho}$ as in \eqref{gamma11rho} and $\eta$,
we have
\begin{equation*}
\begin{split}
\lambda_{\rho} \eta & \sim \max_{T_r\subset \Box: r\geq K_1}\frac{\sharp
\{S: S\subseteq Y^{\rho}_{\Box}\cap T_r\}}{r^{\alpha}}\cdot \sharp \{B:
B\subseteq S\cap Y^{\rho}\} \\
& \lesssim \max_{T_r\subset \Box: r\geq K_1}\frac{\sharp \{B\subset
Y^{\rho}: B\subset T_r\}}{r^{\alpha}} \\
& \leq \frac{\rho^{\frac{m-2}{2}}K^{\frac{3}{2}-\frac{2}{m}}\lambda^{({\rho}
)} (K^{\frac{1}{m}}r)^{\alpha}}{r^{\alpha}},
\end{split}
\end{equation*}
where the last inequality follows from the fact that we can cover a $\rho^{
\frac{m-2}{2}}K^{\frac{1}{2}}r \times K^{\frac{1}{m}}r \times Kr$-tube $T_r$
by $\sim \rho^{\frac{m-2}{2}}K^{\frac{3}{2}-\frac{2}{m}}$ finitely
overlapping $K^{\frac{1}{m}}r$-balls. Hence, we obtain
\begin{equation}  \label{ineq9rho}
\eta \lesssim \frac{\lambda^{(\rho)}\rho^{\frac{m-2}{2}} K^{\frac{3}{2}+
\frac{\alpha-2}{m}}}{\lambda_{\rho}}.
\end{equation}

Finally, we relate $\nu_{\rho }$ and $\nu^{({\rho })}$ by considering the
number of $K$-cubes in each relevant $\frac{R^{1/2}}{\rho^{-\frac{m-2}{2}}}
\times \frac{R^{1/2}}{K^{\frac{1}{2}-\frac{1}{m}}}\times K^{\frac{1}{2}}R^{
\frac{1}{2}}$-tube $S^{\prime }$. Recall that each relevant $S^{\prime }$
contains $\sim \nu _{\rho }$ tubes $S$ in $Y_{\Box }$ and each such $S$
contains $\sim \eta $ $K$-cubes. On the other hand, we can cover $S^{\prime }
$ by $\sim K^{1/2}$ finitely overlapping $R^{\frac{1}{2}}$-cubes, and each $R^{\frac{1}{2}}$-cube contains $\lesssim \nu ^{({\rho })}$
many $K$-cubes in $Y^{\rho }$. Thus, it follows that
\begin{equation}
\nu _{\rho }\lesssim \frac{K^{\frac{1}{2}}\nu ^{({\rho })}}{\eta }.
\label{ineq10rho}
\end{equation}
By inserting \eqref{ineq8rho}, \eqref{ineq9rho} and \eqref{ineq10rho} into
\eqref{decoup7rho}, we derive
\begin{equation*}
\Vert e^{it\phi (D)}f_{\Omega _{\rho }}\Vert _{L^{6}(Y^{\rho })}\lesssim
K^{2\varepsilon ^{2}-\varepsilon }(N^{({\rho })})^{-\frac{1}{3}}(\nu ^{({
\rho })})^{\frac{1}{6}}(\lambda ^{({\rho })})^{\frac{1}{6}}R^{\frac{5}{12}-
\frac{5-\alpha }{6m}+\varepsilon }\Vert f\Vert _{2}.
\end{equation*}
Since $K=R^{\delta }$ and $R$ can be assumed to be sufficiently large
compared to any constant depending on $\varepsilon $, we have
$K^{2\varepsilon ^{2}-\varepsilon }\ll 1$, and the induction closes. Recall that
$\frac{N}{\log K(\log R)^{4}}\lesssim N^{(\rho )}\leq N,\nu ^{(\rho )}\leq
\nu $ and $\lambda ^{(\rho )}\leq \lambda $. This yields
\begin{equation*}
\Vert e^{it\phi (D)}f_{\Omega _{\rho }}\Vert _{L^{6}(Y^{\rho })}\lesssim
_{\varepsilon }N^{-\frac{1}{3}}\lambda ^{\frac{1}{6}}\nu ^{\frac{1}{6}}R^{
\frac{5}{12}-\frac{5-\alpha }{6m}+\varepsilon }\Vert f\Vert _{2}.
\end{equation*}
Combining it with \eqref{case1YandYrho}, we get
\begin{equation*}
\Vert e^{it\phi (D)}f_{\Omega _{1}}\Vert _{L^{6}(Y^{1})}\lesssim
_{\varepsilon }N^{-\frac{1}{3}}\lambda ^{\frac{1}{6}}\nu ^{\frac{1}{6}}R^{
\frac{5}{12}-\frac{5-\alpha }{6m}+\varepsilon }\Vert f\Vert _{2}.
\end{equation*}

\section{\textbf{The Proof of Proposition \ref{mainlemma}}}

To prove Proposition \ref{mainlemma}, we divide $[0,1]^2$ into two parts
\begin{equation*}
[0,1]^2=\Omega_0\bigcup \Omega_1,
\end{equation*}
where $\Omega_0:=[0,1]\times [K^{-\frac{1}{m}},1]$, $\Omega_1:=[0,1]\times
[0,K^{-\frac{1}{m}}]$. We denote $\hat{f}\mid_{\Omega_0}$ and $\hat{f}
\mid_{\Omega_1}$ by $\hat{f}_{\Omega_0}$ and $\hat{f}_{\Omega_1}$,
respectively. Given a $K$-cube $B$, by the triangle inequality,
we have either
$$\Vert e^{it \psi(D)}f_{\Omega_0}
\Vert_{L^6(B)}\geq \frac{1}{8}\Vert e^{it \psi(D)}f\Vert_{L^6(B)}$$
or
$$\Vert e^{it \psi(D)}f_{\Omega_1} \Vert_{L^6(B)}\geq \frac{1}{8}\Vert e^{it \psi(D)}f\Vert_{L^6(B)}.$$
We sort the $K$-cubes in $Y$
as follows. Denote by
\[\{B \subset Y:\Vert e^{it \psi(D)}f_{\Omega_0}\Vert_{L^6(B)}\geq\frac{1}{8}\Vert e^{it \psi(D)}f\Vert_{L^6(B)}\}\]
and
\[\{B \subset Y: \Vert e^{it \psi(D)}f_{\Omega_1}\Vert_{L^6(B)}\geq \frac{1}{8}\Vert e^{it \psi(D)}f\Vert_{L^6(B)}\}\]
by $Y^0$ and $Y^1$, respectively. Clearly, one has $Y=Y^0\cup Y^1$.
If $\sharp \{B: B \subset Y^0 \} \geq \frac{N}{2}$, we
call it $\Omega_0$-case. If $\sharp \{B: B \subset Y^1\} \geq \frac{N}{2}$,
we call it $\Omega_1$-case.

For the $\Omega_0$-case, we further sort the $K$-cubes in $Y^0$ as follows:

\begin{enumerate}
[(1)]
\item For any dyadic number $A^{(0)}$, let $Y^0_{A^{(0)}}$ be the union of
the $K$-cubes in $Y^0$ satisfying
\begin{equation*}
\Vert e^{it\psi(D)}f_{\Omega_0} \Vert_{L^6(B)}\approx A^{(0)}.
\end{equation*}

\item Fix $A^{(0)}$, for any dyadic number $\nu^{(0)}$, let $
Y^0_{A^{(0)},\nu^{(0)}}$ be the union of the $K$-cubes in $Y^0_{A^{(0)}}$
such that for each $B\subset Y^0_{A^{(0)},\nu^{(0)}}$, the $R^{1/2}$-cube
intersecting $B$ contains $\approx \nu^{(0)}$ cubes from $Y^0_{A^{(0)}}$.
\end{enumerate}

With the assumption $\Vert f \Vert_2=1$ the dyadic numbers $A^{(0)},
\nu^{(0)}$ making significant contributions are between $R^{-C}$ and $R^C$.
Therefore, there exist some dyadic numbers $A^{(0)}, \nu^{(0)}$ such that $
\sharp \{B: B \subset Y^0_{A^{(0)},\nu^{(0)}}\} =:N^{(0)} \gtrsim \frac{N}{
(\log R)^2}$. Fix a choice of $A^{(0)}, \nu^{(0)}$ and denote $
Y^0_{A^{(0)},\nu^{(0)}}$ by $Y^0$ for convenience. Then, in the $\Omega_0$
-case we have
\begin{equation*}
\Vert e^{it \psi(D)}f \Vert_{L^6(Y)}\lesssim_{\varepsilon}R^{\varepsilon}
\Vert e^{it \psi(D)}f_{\Omega_0} \Vert_{L^6(Y^0)},
\end{equation*}
and
\begin{equation*}
\frac{N}{(\log R)^2} \lesssim N^{(0)} \leq N, \quad \nu^{(0)} \leq \nu,\quad
\lambda^{(0)} \leq \lambda,
\end{equation*}
where
\begin{equation*}
\lambda^{(0)}:=\max_{B^3(x^{\prime },r)\subset Q_R,x^{\prime }\in \mathbb{R}
^3,r\geq 1}\frac{\sharp\{B_k \subset Y^0: B_k \subset B^3(x^{\prime },r)\}}{
r^{\alpha}}.
\end{equation*}
By a similar argument as in the proof of \eqref{omega0}, we have
\begin{equation*}  \label{omega0l}
\Vert e^{it \psi(D)}f_{\Omega_0} \Vert_{L^6(Y^0)}\leq
C_{\varepsilon}K^{O(1)}(N^{(0)})^{-\frac{1}{3}}(\lambda^{(0)})^{\frac{1}{6}
}(\nu^{(0)})^{\frac{1}{6}}R^{\frac{\alpha}{12}+\varepsilon}\Vert f \Vert_2.
\end{equation*}

For the $\Omega_1$-case, we sort the $K$-cubes in $Y^1$ as similar as in the
$\Omega_0$-case:

\begin{enumerate}
[(1)]
\item For any dyadic number $A^{(1)}$, let $Y^1_{A^{(1)}}$ be the union of
the $K$-cubes in $Y^1$ satisfying
\begin{equation*}
\Vert e^{it \psi(D)}f_{\Omega_1} \Vert_{L^6(B)}\approx A^{(1)}.
\end{equation*}

\item Fix $A^{(1)}$, for any dyadic number $\nu^{(1)}$, let $
Y^1_{A^{(1)},\nu^{(1)}}$ be the union of the $K$-cubes in $Y^1_{A^{(1)}}$
such that for each $B\subset Y^1_{A^{(1)},\nu^{(1)}}$, the $R^{1/2}$-cube
intersecting $B$ contains $\approx \nu^{(1)}$ cubes from $Y^1_{A^{(1)}}$.
\end{enumerate}

The dyadic numbers $A^{(1)}, \nu^{(1)}$ making significant contributions can
be assumed to be between $R^{-C}$ and $R^C$. Therefore, there exist some
dyadic numbers $A^{(1)}, \nu^{(1)}$ such that $\sharp \{B: B \subset
Y^1_{A^{(1)},\nu^{(1)}}\} =:N^{(1)} \gtrsim \frac{N}{(\log R)^2}$ many cubes
$B$. Fix that choice of $A^{(1)}, \nu^{(1)}$ and denote $Y^1_{A^{(1)},
\nu^{(1)}}$ by $Y^1$ for convenience. Then, in the $\Omega_1$-case we have
\begin{equation*}
\Vert e^{it \psi(D)}f \Vert_{L^6(Y)}\lesssim_{\varepsilon}R^{\varepsilon}
\Vert e^{it \psi(D)}f_{\Omega_1} \Vert_{L^6(Y^1)},
\end{equation*}
and
\begin{equation*}
\frac{N}{(\log R)^2} \lesssim N^{(1)} \leq N,\quad \nu^{(1)} \leq \nu,\quad
\lambda^{(1)} \leq \lambda,
\end{equation*}
where
\begin{equation*}
\lambda^{(1)}:=\max_{B^3(x^{\prime },r)\subset Q_R,x^{\prime }\in \mathbb{R}
^3,r \geq K}\frac{\sharp\{B_k \subset Y^1: B_k \subset B^3(x^{\prime },r)\}}{
r^{\alpha}}.
\end{equation*}
We are going to estimate $\Vert e^{it\psi(D)}f_{\Omega_1} \Vert_{L^6(Y^1)}$.
To this end, we partition $\Omega_1$ into disjoint rectangles $\tau$ of
dimensions $K^{-\frac{1}{2}}\times K^{-\frac{1}{m}}$, and write $
f_{\Omega_1}=\sum_{\tau}f_{\tau}$. For each $K$-cube $B$, we have the
following decoupling inequality due to Bourgain-Demeter \cite{BD15} (see
also Lemma 3.3 in \cite{LiZheng21}).

\begin{lemma}
\label{decouplingl} Suppose that $B$ is a $K$-cube in $\mathbb{R}^3$. Then
for any $\varepsilon >0$, there exists a constant $C_{\varepsilon}$ such
that
\begin{equation*}  \label{decoupineql}
\Vert e^{it\psi(D)}f_{\Omega_1} \Vert_{L^6(B)}\leq
C_{\varepsilon}K^{\varepsilon}\Big(\sum_{\tau} \Vert e^{it\psi(D)}f_{\tau}
\Vert^2_{L^6(\omega_{B})} \Big)^{1/2},
\end{equation*}
where $\omega_{B}$ is a weight function essentially supported on $B$.
\end{lemma}

We decompose the function $f_{\Omega_1}$ in physical space further and
perform some dyadic pigeonholing like in Section 2.

First we divide the physical square $[0,R]^2$ into $\frac{R}{K^{1/2}}\times
\frac{R}{K^{1-\frac{1}{m}}}$-rectangles $D$. For each pair $(\tau,D)$, let $
f_{\Box_{\tau,D}}$ be the function formed by cutting off $f$ on the
rectangle $D$ (with a Schwartz tail) in physical space and the rectangle $
\tau$ in Fourier space. We see that $e^{it \psi(D)}f_{\Box_{\tau,D}}$ is
essentially supported on an $\frac{R}{K^{1/2}}\times \frac{R}{K^{1-\frac{1}{m
}}}\times R$-box, which is denoted by $\Box_{\tau,D}$.

\begin{remark}
\label{rek:section3box} In $x_2$ the wave packets follow the evolution of
the phase $\xi^m_2$, but in $x_1$ they essentially follow the
Schr\"{o}dinger evolution because the phase $\phi_1(\xi_1)$ is basically like
$\xi^2_1 $. So the scaling in both directions is different.
\end{remark}

The box $\Box_{\tau,D}$ is parallel to the normal direction of the surface $
\Sigma_{\tau}$ at the left bottom corner of $\tau$. For any fixed $\tau$,
the different boxes $\Box_{\tau,D}$ tile $Q_R$. In particular, for each $
\tau $, a given $K$-cube $B$ lies in exactly one box $\Box_{\tau,D}$. We
have
\begin{equation*}
f_{\Omega_1}=\sum_{\tau}\sum_{D}f_{\Box_{\tau,D}}
\end{equation*}
and write $f_{\Omega_1}=\sum_{\Box}f_{\Box}$ for abbreviation. By Lemma \ref
{decouplingl}, for each $K$-cube $B$, there holds
\begin{equation*}  \label{decoupineq1l}
\Vert e^{it\psi(D)}f_{\Omega_1} \Vert_{L^6(B)}\leq
C_{\varepsilon}K^{\varepsilon^2}\Big(\sum_{\Box} \Vert e^{it\psi(D)}f_{\Box}
\Vert^2_{L^6(\omega_{B})} \Big)^{1/2}.
\end{equation*}
Recall that $K=R^{\delta}$, where $\delta=\varepsilon^{100}$. We denote by
\begin{equation*}
R_1:=\frac{R}{K}=R^{1-\delta}, \quad K_1=R^{\delta}_1=R^{\delta-\delta^2}.
\end{equation*}
Tile $\Box$ by $K^{\frac{1}{2}}K_1\times K^{\frac{1}{m}}K_1 \times K K_1$
-tube $S$, and also tile $\Box$ by $R^{\frac{1}{2}}\times \frac{R^{1/2}}{K^{
\frac{1}{2}-\frac{1}{m}}}\times K^{\frac{1}{2}}R^{\frac{1}{2}}$-tubes $
S^{\prime }$ (all running parallel to the long axis of $\Box$). After rescaling
the $\Box$ becomes an $R_1$-cube, the tubes $S^{\prime }$ and $S$ become
lattice $R_1^{\frac{1}{2}}$-cubes and $K_1$-cubes, respectively. We regroup
tubes $S$ and $S^{\prime }$ inside each $\Box$ as follows:

\begin{enumerate}
[(1)]
\item Sort those tubes $S$ which intersect $Y^1$ according to the value $
\Vert e^{it\psi(D)}f_{\Box}\Vert_{L^6(S)}$ and the number of $K$-cubes
contained in it. For dyadic numbers $\eta,\beta_1$, we use $\mathbb{S}
_{\Box,\eta,\beta_1}$ to stand for the collection of tubes $S\subset \Box$
each of which containing $\sim \eta$ $K$-cubes and $\Vert
e^{it\psi(D)}f_{\Box}\Vert_{L^6(S)}\sim \beta_1$.

\item For fixed $\eta,\beta_1$, we sort the tubes $S^{\prime }\subset \Box$
according to the number of the tubes $S\in \mathbb{S}_{\Box,\eta,\beta_1}$
contained in it. For dyadic number $\nu_1$, let $\mathbb{S}
_{\Box,\eta,\beta_1,\nu_1}$ be the subcollection of $\mathbb{S}
_{\Box,\eta,\beta_1}$ such that for each $S\in \mathbb{S}_{\Box,\eta,
\beta_1,\nu_1}$, the tube $S^{\prime }$ containing $S$ contains $\sim \nu_1$
tubes from $\mathbb{S}_{\Box,\eta,\beta_1}$.

\item For fixed $\eta,\beta_1,\nu_1$, we sort the boxes $\Box$ according to
the value $\Vert f_{\Box} \Vert_2$, the number $\sharp \mathbb{S}
_{\Box,\eta,\beta_1,\nu_1}$ and the value $\lambda_1$ defined below. For
dyadic numbers $\beta_2,N_1,\lambda_1$, let $\mathbb{B}_{\eta,\beta_1,\nu_1,
\beta_2,N_1,\lambda_1}$ denote the collection of boxes $\Box$ each of which
satisfying that
\begin{equation*}
\Vert f_{\Box} \Vert_2\sim \beta_2,\quad \sharp \mathbb{S}
_{\Box,\eta,\beta_1,\nu_1}\sim N_1
\end{equation*}
and
\begin{equation}  \label{gamma11l}
\max_{T_r\subset \Box:r\geq K_1}\frac{\sharp\{S\in \mathbb{S}
_{\Box,\eta,\beta_1,\nu_1}:S\subset T_r\}}{r^{\alpha}}\sim\lambda_1,
\end{equation}
where $T_r$ are $K^{\frac{1}{2}}r\times K^{\frac{1}{m}}r \times Kr$-tubes in
$\Box$ running parallel to the long axis of $\Box$.
\end{enumerate}

Let $Y^1_{\Box,\eta,\beta_1,\nu_1}$ denote $\{S:S\subset \mathbb{S}
_{\Box,\eta,\beta_1,\nu_1}\}$. By pigeonholing, we can choose
$\eta,\beta_1,\nu_1,\beta_2,N_1,\lambda_1$ such that
\begin{equation}  \label{decoup314l}
\begin{split}
& \quad \Vert e^{it \psi(D)}f_{\Omega_1} \Vert_{L^6(B)} \\
& \lesssim (\log R)^6 K^{\varepsilon^4}\Big(\sum_{\Box \in \mathbb{B}
_{\eta,\beta_1,\nu_1,\beta_2,N_1,\lambda_1},B\subset
Y_{\Box,\eta,\beta_1,\nu_1}}\Vert e^{it\psi(D)}f_{\Box}
\Vert^2_{L^6(\omega_{B})} \Big)^{1/2}
\end{split}
\end{equation}
holds for a fraction $\gtrsim (\log R)^{-6}$ of all $K$-cubes $B$. For
brevity, we denote by
\begin{equation*}
Y^1_{\Box}:=Y^1_{\Box,\eta,\beta_1,\nu_1}, \quad \mathbb{B}:=\mathbb{B}
_{\eta,\beta_1,\nu_1,\beta_2,N_1,\lambda_1}.
\end{equation*}
Finally, we sort the $K$-cubes $B$ satisfying \eqref{decoup314l} by $\sharp
\{\Box\in \mathbb{B}: B\subset Y^1_{\Box}\}$. Let $Y^{\prime }_{\mu}\subset
Y^1$ be a union of cubes that obey
\begin{equation}  \label{decoup315l}
\Vert e^{it\psi(D)}f_{\Omega_1} \Vert_{L^6(B)} \lesssim (\log R)^6
K^{\varepsilon^2}\Big(\sum_{\Box \in \mathbb{B},B\subset Y^1_{\Box}}\Vert
e^{it \psi(D)}f_{\Box} \Vert^2_{L^6(\omega_{B})} \Big)^{1/2}
\end{equation}
and
\begin{equation}  \label{boxdyadicl}
\sharp \{\Box\in \mathbb{B}: B\subset Y^1_{\Box}\}\sim \mu.
\end{equation}
Since the dyadic number $\mu$ must satisfy $1\leq \mu \leq R^{C}$, there are
$\sim\log R$ different choices for $\mu$, so by pigeonholing there exists
one such $\mu$ such that the number of cubes in $Y^{\prime }_{\mu}$ is $
\gtrsim \frac{N}{(\log R)^{7}}$. We have
\begin{equation}  \label{ineq6l}
\Vert e^{it\psi(D)}f_{\Omega_1} \Vert^6_{L^6(Y)} \lesssim (\log R)^7
\sum_{B\subset Y^{\prime }} \Vert e^{it\psi(D)}f_{\Omega_1} \Vert^6_{L^6(B)}.
\end{equation}
For each $B\subset Y^{\prime }$, it follows from \eqref{decoup315l},
\eqref{boxdyadicl} and H\"{o}lder's inequality that
\begin{equation}  \label{decoup316l}
\Vert e^{it\psi(D)}f_{\Omega_1} \Vert^6_{L^6(B)} \lesssim (\log R)^{36}
K^{6\varepsilon^2}\mu^2 \sum_{\Box \in \mathbb{B},B\subset Y_{\Box}}\Vert
e^{it\psi(D)}f_{\Box} \Vert^6_{L^6(\omega_{B})}.
\end{equation}
Putting \eqref{ineq6l} and \eqref{decoup316l} together, one has
\begin{equation}  \label{decoup7l}
\Vert e^{it\psi(D)}f_{\Omega_1} \Vert_{L^6(Y^1)} \lesssim (\log R)^{C}
K^{\varepsilon^2}\mu^{1/3}\Big(\sum_{\Box \in \mathbb{B}}\Vert e^{it
\psi(D)}f_{\Box} \Vert^6_{L^6(Y^1_{\Box})} \Big)^{1/6}.
\end{equation}
The key point is that phase functions of the form $\psi(\xi_1,\xi_2)$ are
closed over each region $\tau= \tau_{\Box}$ in $\Omega_1$ under the change
of variable
\begin{equation*}
\left\{
\begin{array}{l}
\xi_1 = a + K^{-\frac{1}{2}}\eta_1,\;a \in [\frac{1}{2}, 1-K^{-\frac{1}{2}
}]\cap K^{-1/2}\mathbb{Z}, \\
\xi_2 = K^{-\frac{1}{m}}\eta_2.
\end{array}
\right.
\end{equation*}
More precisely, under the change of variables $\psi(\xi_1,\xi_2)$ becomes
\begin{equation*}
\tilde{\psi}(\eta_1,\eta_2):=K[\phi_1(a+K^{-1/2}\eta_1)-\phi_1(a)-K^{-1/2}
\eta_1 \phi^{\prime }_1(a)]+\eta^m_2,
\end{equation*}
which satisfies the condition \eqref{condition}. We write
\begin{equation*}
\vert e^{it\psi(D)}f_{\Box}(x)\vert= K^{-\frac{m+2}{4m}}\vert e^{i\tilde{t}
\tilde{\psi}(D)}g(\tilde{x})\vert,
\end{equation*}
where
\begin{equation*}
\hat{g}(\eta_1,\eta_2):= K^{-\frac{m+2}{4m}}\hat{f}(a+K^{-\frac{1}{2}
}\eta_1,K^{-\frac{1}{m}}\eta_2),\quad \Vert g \Vert_2= \Vert f \Vert_2,
\end{equation*}
and
\begin{equation*}
\left\{
\begin{array}{l}
\tilde{x}_1 := K^{-\frac{1}{2}}(x_1+\phi^{\prime }_1(a)t), \\
\tilde{x}_2 := K^{-\frac{1}{m}}x_2, \\
\tilde{t} := K^{-1}t.
\end{array}
\right.
\end{equation*}
We denote the above relation by $(\tilde{x},\tilde{t})=\mathcal{L}(x,t)$ and
write
\begin{equation}  \label{relationl}
\Vert e^{it\psi(D)}f_{\Box} \Vert_{L^6(Y^1_{\Box})}=K^{-\frac{1}{3m}}\Vert
e^{i\tilde{t}\tilde{\psi}(D)}g(\tilde{x}) \Vert_{L^6(\tilde{Y}^1)},
\end{equation}
where $\tilde{Y}^1=\mathcal{L}(Y^1_{\Box})$. By \eqref{relationl} and the
inductive hypothesis at scale $R_1$, we derive
\begin{equation*}  \label{inductionl}
\Vert e^{it\psi(D)}f_{\Box} \Vert_{L^6(Y^1_{\Box})}\lesssim K^{-\frac{1}{3m}
}N_1^{-\frac{1}{3}}\lambda_1^{\frac{1}{6}}\nu_1^{\frac{1}{6}}(\tfrac{R}{K})^{
\frac{1}{3}-\frac{4-\alpha}{6m}+\varepsilon}\Vert f_{\Box} \Vert_2.
\end{equation*}

Now we consider the cardinality of the set $\{(\Box,B): \Box \in \mathbb{B},
B\subseteq Y^1_{\Box}\cap Y^{\prime }\}$. On one hand, by the choice of $\mu$
as in \eqref{boxdyadicl}, there is a lower bound
\begin{equation*}
\sharp \{(\Box,B): \Box \in \mathbb{B}, B\subseteq Y^1_{\Box}\cap Y^{\prime
-7}N^{(1)}\mu.
\end{equation*}
On the other hand, by the choices of $N_1$ and $\eta$, for each $\Box \in
\mathbb{B}$, $Y^1_{\Box}$ contains $\sim N_1$ tubes $S$ and each $S$
contains $\sim \eta$ cubes in $Y^1$, so
\begin{equation*}
\sharp \{(\Box,B): \Box \in \mathbb{B}, B\subseteq Y^1_{\Box}\cap Y^{\prime
}\}\lesssim (\sharp \mathbb{B})N_1\eta.
\end{equation*}
Therefore, we get
\begin{equation}  \label{ineq8l}
\frac{\mu}{\sharp \mathbb{B}} \lesssim \frac{(\log R)^7 N_1 \eta}{N^{(1)}}.
\end{equation}

Next, by our choices of $\lambda _{1}$ as in \eqref{gamma11l} and $\eta $,
it holds that
\begin{equation*}
\begin{split}
\lambda _{1}\eta & \sim \max_{T_{r}\subset \Box :r\geq K_{1}}\frac{\sharp
\{S:S\subseteq Y_{\Box }^{1}\cap T_{r}\}}{r^{\alpha }}\cdot \sharp
\{B:B\subseteq S\cap Y^{1}\} \\
& \lesssim \max_{T_{r}\subset \Box :r\geq K_{1}}\frac{\sharp \{B\subset
Y:B\subset T_{r}\}}{r^{\alpha }} \\
& \leq \frac{K^{\frac{3}{2}-\frac{2}{m}}\lambda ^{(1)}(K^{\frac{1}{m}
}r)^{\alpha }}{r^{\alpha }},
\end{split}
\end{equation*}
where the last inequality follows from the fact that one can cover a $K^{
\frac{1}{2}}r\times K^{\frac{1}{m}}r\times Kr$-tube $T_{r}$ by $\sim K^{
\frac{3}{2}-\frac{2}{m}}$ finitely overlapping $K^{\frac{1}{m}}r$-balls.
Thus, we obtain
\begin{equation}
\eta \lesssim \frac{\lambda ^{(1)}K^{\frac{3}{2}+\frac{\alpha -2}{m}}}{
\lambda _{1}}.  \label{ineq9l}
\end{equation}

Finally, we relate $\nu _{1}$ and $\nu ^{(1)}$ by considering the number of $
K$-cubes in each relevant $R^{\frac{1}{2}}\times \frac{R^{1/2}}{K^{\frac{1}{2
}-\frac{1}{m}}}\times K^{\frac{1}{2}}R^{\frac{1}{2}}$-tube $S^{\prime }$. On
one hand, each relevant $S^{\prime }$ contains $\sim \nu _{1}$ tubes $S$ in $
Y_{\Box }$ and each such $S$ contains $\sim \eta $ $K$-cubes. On the other
hand, a relevant $S^{\prime }$ can be covered by $\sim K^{1/2}$ finitely
overlapping $R^{\frac{1}{2}}$-cubes and each $R^{\frac{1}{2}}$-cube contains
$\lesssim \nu ^{(1)}$ many $K$-cubes in $Y^{1}$. It follows that
\begin{equation}
\nu _{1}\lesssim \frac{K^{\frac{1}{2}}\nu ^{(1)}}{\eta }.  \label{ineq10l}
\end{equation}
Using the similar argument as in the proof of \eqref{omega3}, we derive
\begin{equation*}
\Vert e^{it\psi (D)}f_{\Omega _{1}}\Vert _{L^{6}(Y^{1})}\lesssim
K^{2\varepsilon ^{2}-\varepsilon }(N^{(1)})^{-\frac{1}{3}}(\nu ^{(1)})^{
\frac{1}{6}}(\lambda ^{(1)})^{\frac{1}{6}}R^{\frac{1}{3}-\frac{4-\alpha }{6m}
+\varepsilon }\Vert f\Vert _{2}.
\end{equation*}
Since $K=R^{\delta }$ and $R$ can be assumed to be sufficiently large
compared to any constant depending on $\varepsilon $, we have $
K^{2\varepsilon ^{2}-\varepsilon }\ll 1$ and the induction closes. Recall
that $\frac{N}{\log R}\lesssim N^{(1)}\leq N,\nu ^{(1)}\leq \nu $ and $
\lambda ^{(1)}\leq \lambda $. This yields
\begin{equation*}
\Vert e^{it\psi (D)}f_{\Omega _{1}}\Vert _{L^{6}(Y^{1})}\lesssim
_{\varepsilon }N^{-\frac{1}{3}}\lambda ^{\frac{1}{6}}\nu ^{\frac{1}{6}}R^{
\frac{1}{3}-\frac{4-\alpha }{6m}+\varepsilon }\Vert f\Vert _{2}.
\end{equation*}
We complete the proof of Lemma \ref{mainlemma}.

\section{\textbf{Applications of Corollary \ref{weakfractal}}}

\subsection{Application to the average Fourier decay of fractal measures associated with the surfaces $F^2_m$.}
We recall the definition of $\alpha$-dimensional measure as follows.

\begin{definition}
\label{measure} Let $\alpha \in (0,d]$. We say that $\mu$ is an $\alpha$
-dimensional measure in $\mathbb{R}^d$ if it is a probability measure
supported in the unit ball $B^{d}(0,1)$ and satisfies
\begin{equation*}
\mu(B(x,r))\leq C_{\mu}r^{\alpha},\;\forall r>0,\;\forall x \in \mathbb{R}^d.
\end{equation*}
\end{definition}

We denote $d\mu _{R}(\cdot ):=R^{\alpha }d\mu (\frac{\cdot }{R})$. Let $
\gamma _{d}(\alpha )$ denote the supremum of the numbers $\gamma $ for which
\begin{equation*}
\Vert \hat{\mu}(R\cdot )\Vert _{L^{2}(\mathbb{S}^{d-1})}^{2}\leq C_{\alpha
,\mu }R^{-\gamma }
\end{equation*}
whenever $R>1$ and $\mu $ is an $\alpha $-dimensional measure in $\mathbb{R}
^{d}$. The problem of identifying the precise value of $\gamma _{d}(\alpha )$
was proposed by Mattila \cite{Mattila04}. The lower bound of $\gamma
_{d}(\alpha )$ has been studied by several authors. For instance, one can
see ~\cite{Wolff1999, DGOWWZ18, DuZhang19-Annals}.

Let $\beta_{2,m}(\alpha )$ denote the supremum of the numbers $\beta $ for
which
\begin{equation*}
\Vert \hat{\mu}(R\cdot )\Vert _{L^{2}(F_{m}^{2})}^{2}\leq C_{\alpha ,\mu
}R^{-\beta }
\end{equation*}
whenever $R>1$ and $\mu $ is an $\alpha $-dimensional measure in $\mathbb{R}
^{3}$. We will study the lower bound of $\beta _{2,m}(\alpha )$. For $m=2$,
the surface $F_{m}^{2}$ is exactly the paraboloid $P^{2}$ over the region $
[0,1]^{2}$. Note that the Gaussian curvature of the surface $F_{m}^{2}$
vanishes when $\xi _{1}=0$ or $\xi _{2}=0$, which is different from the
sphere or the paraboloid case studied in the literature.

We are going to derive the following result.

\begin{theorem}
\label{mainthmdecay} Let $m\geq 4$ be an even number and $0< \alpha \leq 3$.
Then
\begin{equation*}  \label{maindecayineq}
\beta_{2,m}(\alpha)\geq \big(\frac{5}{6}-\frac{1}{3m}\big)\alpha+\frac{5}{3m}-
\frac{5}{6}.
\end{equation*}
\end{theorem}

By the arguments in Remark 2.5 from \cite{DGOWWZ18}, Theorem \ref
{mainthmdecay} can be reduced to the weighted restriction estimate in
Theorem \ref{weightedrestrictionthm}.
In fact, one can obtain
$$\beta_{2,m}(\alpha)\geq 2\big(\frac{\alpha}{2}-\gamma\big),$$
if  there holds that  for any $\varepsilon>0,$
$$
\Vert  e^{it\phi(D)}f \Vert_{L^2(B^3(0,R);d\mu_{R}(x,t))}\leq C_{\varepsilon}  R^{\gamma+\varepsilon }\Vert f \Vert_2.
$$

\begin{theorem}
\label{weightedrestrictionthm} Let $m\geq 4$ be an even number and $0<\alpha
\leq 3$. Suppose that $\mu $ is an $\alpha $-dimensional measure in $\mathbb{
R}^{3}$. For all $R>1$ and any $\varepsilon >0$, there exists a positive
constant $C_{\varepsilon }$ such that
\begin{equation*}
\Vert e^{it\phi (D)}f\Vert _{L^{2}(B^{3}(0,R);d\mu _{R}(x,t))}\leq
C_{\varepsilon }R^{\frac{5}{12}-\frac{5}{6m}+(\frac{1}{12}+\frac{1}{6m}
)\alpha +\varepsilon }\Vert f\Vert _{2},
\end{equation*}
holds for all $f$ with Fourier supports in $[0,1]^{2}$.
\end{theorem}

Using Theorem \ref{fractalrestriction} and the argument in the proof of
Theorem 2.2 from \cite{DuZhang19-Annals} we may prove Theorem \ref
{weightedrestrictionthm} as follows.

Denote $e^{it\phi (D)}f(x)$ by $\mathcal{E}f(x,t)$, and $(x,t)$ by $\tilde{x}
$. Since $supp\;\hat{f}\subseteq B^{2}(0,1)$, we have $supp\;\widehat{
\mathcal{E}f}\subseteq B^{3}(0,1)$. Thus there exists a Schwartz bump
function $\psi $ on $\mathbb{R}^{3}$ such that $(\mathcal{E}f)^{2}=(\mathcal{
E}f)^{2}\ast \psi $. The function $\max_{|\tilde{y}-\tilde{x}|\leq
e^{300}}|\psi (\tilde{y})|,$ which we denote it by $\psi _{1}(\tilde{x}),$
rapidly decays. Note also that any $(x,t)$ in $\mathbb{R}^{3}$ belongs to a
unique integer lattice cube whose center we denote by $\tilde{m}
=(m,m_{3})=(m_{1},m_{2},m_{3})=\tilde{m}(x,t).$ Then we have
\begin{equation*}
\begin{split}
& \quad \left\Vert e^{it\phi (D)}f\right\Vert _{L^{2}(B^{3}(0,R);d\mu
_{R}(x,t))}^{2} \\
& =\int_{B^{3}(0,R)}|\mathcal{E}f(x,t)|^{2}d\mu _{R}(x,t) \\
& \leq \int_{B^{3}(0,R)}(|\mathcal{E}f|^{2}\ast \psi )(x,t)d\mu _{R}(x,t) \\
& \lesssim \int_{B^{3}(0,R)}(|\mathcal{E}f|^{2}\ast \psi _{1})(\tilde{m}
(x,t))d\mu _{R}(x,t) \\
& \leq \sum_{\tilde{m}\in \mathbb{Z}^{3},|m_{i}|\leq R}\Big(\int_{|(x,t)-
\tilde{m}|\leq 10}d\mu _{R}(x,t)\Big)\cdot (|\mathcal{E}f|^{2}\ast \psi
_{1})(\tilde{m}).
\end{split}
\end{equation*}
Without loss of generality, we may assume $\Vert f\Vert _{2}=1$. For each $
\tilde{m}\in \mathbb{Z}^{3}$ we define
\begin{equation*}
\nu _{\tilde{m}}:=\int_{|(x,t)-\tilde{m}|\leq 10}d\mu _{R}(x)\lesssim 1.
\end{equation*}
From the discussion above, we have
$(\vert\mathcal{E}f\vert^2 \ast \psi_1)(\tilde{m})\lesssim 1, $
which implies
\begin{equation*}
\Vert e^{it\phi (D)}f\Vert _{L^{2}(B^{3}(0,R);d\mu _{R}(x,t))}^{2}\lesssim
\sum_{\nu \in \lbrack R^{-300},1]}\sum_{\tilde{m}:\nu _{\tilde{m}\sim \nu
}}\nu \cdot (|\mathcal{E}f|^{2}\ast \psi _{1})(\tilde{m})+R^{-270},
\end{equation*}
where the first summation is over all dyadic numbers $\nu $. For each dyadic
$\nu $, denote $A_{\nu }=\{\tilde{m}\in \mathbb{Z}^{3}:\;|m_{i}|\leq R,\nu _{
\tilde{m}}\sim \nu \}$. Performing a dyadic pigeonholing over $\nu $ we see
that there exists a dyadic $\nu \in \lbrack R^{-300},1]$ such that for any
small $\varepsilon >0$,
\begin{equation*}
\begin{split}
& \quad \left\Vert e^{it\phi (D)}f\right\Vert _{L^{2}(B^{3}(0,R);d\mu
_{R}(x,t))}^{2} \\
& \lesssim _{\varepsilon }R^{\varepsilon }\sum_{\tilde{m}\in A_{\nu }}\nu
\left( \int_{B^{3}(\tilde{m},R^{\varepsilon })}|\mathcal{E}f|^{2}\right)
+R^{-240} \\
& \lesssim _{\varepsilon }R^{\varepsilon }\nu \cdot \int_{\bigcup_{\tilde{m}
\in A_{\nu }}B^{3}(\tilde{m},R^{\varepsilon })}|\mathcal{E}f|^{2}+R^{-240}.
\end{split}
\end{equation*}
Consider the set $X_{\nu}=\bigcup_{m \in A_{\nu}}B^3(\tilde{m},R^{\varepsilon})$. For every $r> R^{2\varepsilon}$, by the definition of $A_{\nu}$ and the assumption of $\mu_R$, for any $\tilde x\in B^3(0,r)$, we have
\[
\nu \sum_{\tilde m\in A_\nu
\atop B^3(\tilde m,R^\varepsilon) \subset B^3(\tilde x,r)}1 \leq  \sum_{\tilde m\in A_\nu}
\int_{B^3(\tilde m,R^\varepsilon) \cap B^3(\tilde x,r)} d\mu_R \leq \int_{B^3(\tilde x,r)} d\mu_R \leq r^\alpha.
\]
Thus,
the intersection of $X_{\nu}$ and any $r-$ball can be contained in no more than $\nu^{-1}r^{\alpha}$ disjoint $R^{\varepsilon}$-balls.
Hence, we can decompose  the set $
\bigcup_{\tilde{m}\in A_{\nu}}B^3(\tilde{m},R^{\varepsilon})$ further into a union of unit cubes $X=\cup B_k$ satisfying that
\[
\#\{ B_k\subset X: B_k\subset B^3(\tilde x,r) \} \lesssim R^{30 \epsilon} \nu ^{-1}r^\alpha\text{  for any } \tilde x\in \mathbb R^3\text{ and } r\geq 1.\]
Now,  we can apply Corollary \ref{weakfractal} to $X$ with $\lambda \lesssim R^{30\varepsilon}\nu^{-1}$ and $\alpha$, which yields
\begin{equation*}
  \begin{split}
  \quad  \left\Vert e^{it\phi(D)}f \right\Vert_{L^2(B^3(0,R);d\mu_R(x,t))} & \lesssim \nu^\frac12 \left\Vert e^{it\phi(D)}f \right\Vert_{L^2(X)} \\
    & \lesssim  C_{\varepsilon}\nu^{\frac12-\frac{1}{3}}R^{\frac{5}{6}+\frac{\alpha}{6}-\frac{5-\alpha}{3m}+20\varepsilon}\Vert f \Vert_{L^2}\\
  & \lesssim C_{\varepsilon}R^{\frac{5}{6}+\frac{\alpha}{6}-\frac{5-\alpha}{3m}+20\varepsilon}\Vert f \Vert_{L^2}.
   \end{split}
\end{equation*}
Here, we used the fact $\nu\lesssim 1$ in the last inequality.

It completes the proof of Theorem \ref{weightedrestrictionthm}.

\subsection{Application to the almost everywhere convergence problem  associated to fractal  dimensional measures}
One can also consider the best Sobolev exponent for the almost everywhere
convergence problems associated to $\alpha $-dimensional Hausdorff measures.
Results of this direction are widely studied by many authors, and one can refer
to \cite{BBCR-2011,DuGuthLiZhang,DuZhang19-Annals,Mattila-2015,EV2021} for
more details.

We denote
\begin{equation*}
s_{m}(\alpha )=\inf \Big\{s\geq 0:\lim_{t\rightarrow 0}e^{it\phi
(D)}f=f~\alpha \text{-a.e.}\text{ for every }f\in H^{s}(\mathbb{R}^{2})\Big\}.
\end{equation*}
Combining Corollary \ref{a1estimate} and results of Eceizabarrena and
Ponce-Vanegas \cite{EV2021} and the trivial dispersive estimates
\begin{equation*}
\Vert e^{it\phi (D)}\varphi (x)\Vert _{L_{x}^{\infty }(\mathbb{R}^{2})}\leq
C_{\varphi }~|t|^{-\frac{1}{2}-\frac{1}{m}}
\end{equation*}
for any Schwartz function $\varphi $ with support in $A(1)$ \footnote{
In fact, we need Theorem 1.1, Theorem 1.2, and Lemma 2.4 in \cite{EV2021}.},
we have

\begin{corollary}
\label{fractal-divergence-2} For $m\geq 4$ being an even number, we have
\begin{equation*}
s_{m}(\alpha )\left\{ \begin{aligned} =&~\frac{2-\alpha}{2},& 0\leq
\alpha\leq \frac12+\frac1m; \\ \leq &~\frac{3}{4}-\frac{1}{2m},&
\frac12+\frac1m \leq \alpha\leq \frac{7m-2}{5m-2}; \\ \leq &~
\frac{4}{3}-\frac{2}{3m}-\frac{\alpha}{12}\big(5-\frac{2}{m}\big) ,&
\frac{7m-2}{5m-2}\leq \alpha\leq 2. \\ \end{aligned}\right.
\end{equation*}
\end{corollary}

\begin{figure}[H]
\includegraphics[width=12cm]{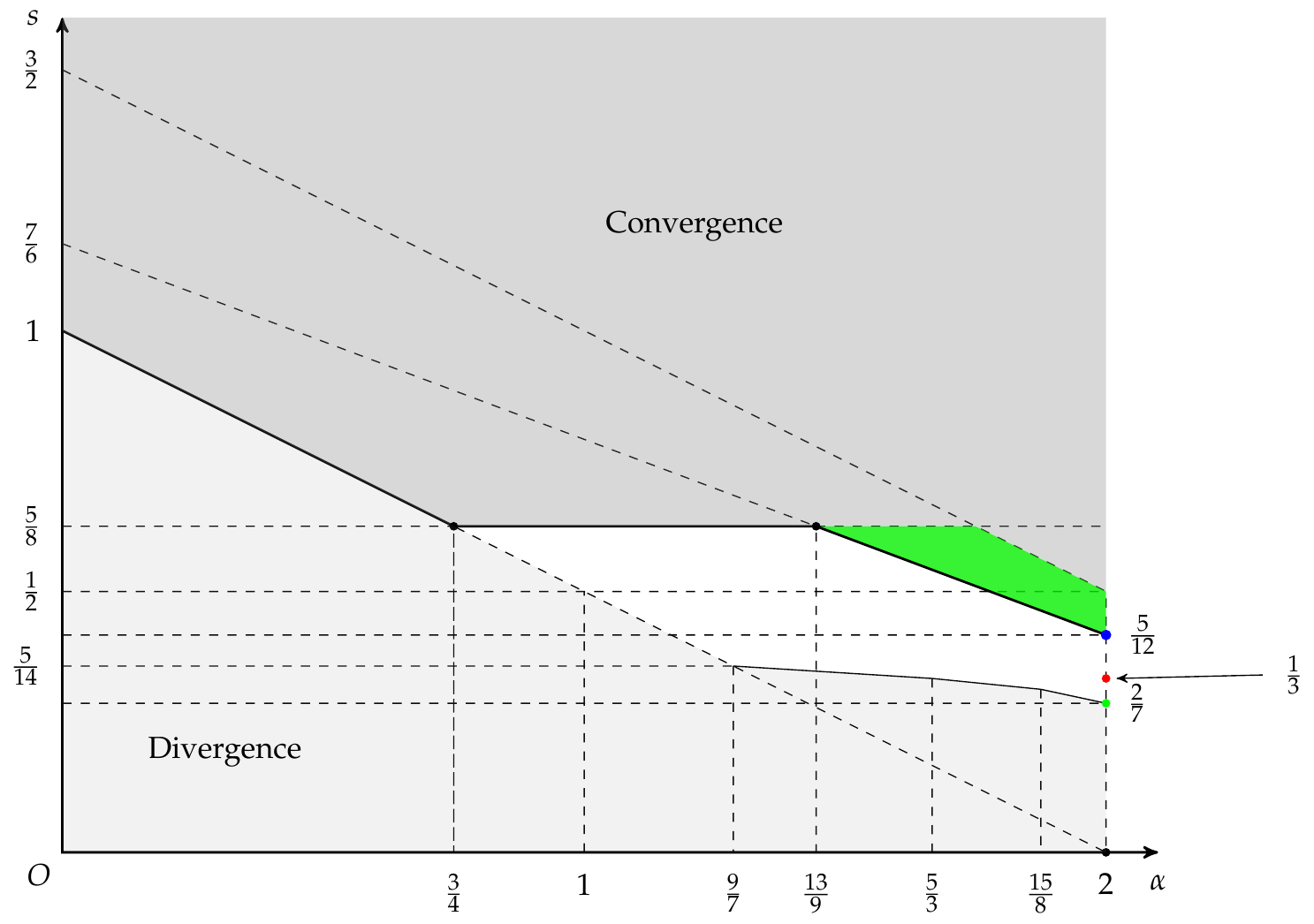}    \caption{Convergence and Divergence regions($m=4$).} \label{figure-m=4}
\end{figure}

\begin{remark}
For general $0<\alpha \leq 2$, we collect the known results of upper bound
and lower bound of $s_{m}(\alpha )$ and show them in Figure \ref{figure-m=4} (we
take $m=4$ for simplicity). Eceizabarrena and Ponce-Vanegas \cite
{EV2021} proved that the light grey region is the divergence region and the
dark grey region is the convergence region. Our results show that the yellow
region is also a convergence region. To the authors' best knowledge, the
convergence of the remaining regions remains open.
\end{remark}

\bigskip \textbf{Acknowledgments.} We thank Daniel Eceizabarrena for
pointing out their results in \cite{EV2021} to us. J. Zhao is supported by
National Natural Science Foundation of China (Grant No. 12101562) and
Natural Science Foundation of Zhejiang (No. LQ20A010003). T. Zhao is
supported by the Fundamental Research Funds for the Central Universities
(FRF-TP-20-076A1).

\bibliographystyle{amsplain}

\end{document}